\documentclass[12pt]{amsart}

\usepackage{am ssymb, stmaryrd, comment}
\usepackage{hyperref}
\usepackage[margin=1.0 in]{geometry}
\usepackage{color,soul}
\usepackage{multicol}
\usepackage[table]{xcolor}
\usepackage{enumerate}
\usepackage{mdframed}
\usepackage[T1]{fontenc}
\usepackage{tikz}
\usetikzlibrary{positioning}
\usepackage{multicol}
\usepackage{caption}
\usepackage{subcaption}
\usepackage{float}

\makeatletter
\def\bbordermatrix#1{\begingroup \m@th
  \@tempdima 4.75\p@
  \setbox\z@\vbox{%
    \def\cr{\crcr\noalign{\kern2\p@\global\let\cr\endline}}%
    \ialign{$##$\hfil\kern2\p@\kern\@tempdima&\thinspace\hfil$##$\hfil
      &&\quad\hfil$##$\hfil\crcr
      \omit\strut\hfil\crcr\noalign{\kern-\baselineskip}%
      #1\crcr\omit\strut\cr}}%
  \setbox\tw@\vbox{\unvcopy\z@\global\setbox\@ne\lastbox}%
  \setbox\tw@\hbox{\unhbox\@ne\unskip\global\setbox\@ne\lastbox}%
  \setbox\tw@\hbox{$\kern\wd\@ne\kern-\@tempdima\left[\kern-\wd\@ne
    \global\setbox\@ne\vbox{\box\@ne\kern2\p@}%
    \vcenter{\kern-\ht\@ne\unvbox\z@\kern-\baselineskip}\,\right]$}%
  \null\;\vbox{\kern\ht\@ne\box\tw@}\endgroup}
\makeatother

\def\VR{\kern-\arraycolsep\strut\vrule &\kern-\arraycolsep}
\def\vr{\kern-\arraycolsep & \kern-\arraycolsep}

\theoremstyle{plain}
\newtheorem{theorem}[subsection]{Theorem}
\newtheorem{prop}[subsection]{Proposition}
\newtheorem{Keylemma}[subsection]{Key Lemma}

\newtheorem{lemma}[subsection]{Lemma}

\newtheorem{cor}[subsection]{Corollary}

\theoremstyle{definition}
\newtheorem{definition}[subsection]{Definition}

\newtheorem{formula}[subsection]{Formula}
\newtheorem{notation}[subsection]{Notation}

\newtheorem{example}[subsection]{Example}

\newtheorem{conjecture}[subsection]{Conjecture}
\newtheorem{remark}[subsection]{Remark}

\normalfont
\usepackage[T1]{fontenc}{}

\newcommand{\m}{\mathfrak{m}}
\newcommand{\PSym}{\mathbb P({\text{Sym}}^d(k^2)^*)}
\newcommand{\PSymN}{\mathbb P({\text{Sym}}^d(k^n)^*)}
\newcommand{\fpt}{\operatorname{fpt}}

\date{\today}

\title{Values of the F-pure threshold for \\ homogeneous polynomials\footnote{Math Subject Classifications: Primary 13A35. Secondary 14B05.}}
\author{ Karen E. Smith and Adela Vraciu}
\email{ kesmith@umich.edu}
\thanks{
This work is partially funded  by NSF DMS \#1952399 and \#2101075.}

\begin{document}

\begin{abstract} We find a formula, in terms of $n, d$ and $p$,
 for the value of the $F$-pure threshold for the generic homogeneous polynomial of degree $d$ in $n$ variables over an algebraically closed field of characteristic $p$. We also show that, in every characteristic $p$ and
  for all $d\geq 4$ not divisible by $p$, there {\it always} exist reduced polynomials of degree $d$ in $k[x, y]$ whose $F$-pure threshold is  a truncation of the base $p$ expansion of $\frac{2}{d}$ at some place; in particular,  there always exist reduced
   polynomials $f$ whose $F$-pure threshold is {\it strictly} less than $\frac{2}{\deg(f)}$.
   We provide an example to resolve, negatively, a question proposed by Hernandez, N\'u\~nez-Betancourt, Witt and Zhang, as to whether a list of necessary restrictions they prove on the $F$-pure threshold of reduced forms are "minimal" for $p\gg 0$.
    On the other hand, we also provide evidence supporting and refining their ideas, including identifying specific truncations of the base $p$ expansion of $\frac{2}{d}$ that are always  $F$-pure thresholds for reduced forms of degree $d$, and computations that 
    show their conditions suffice (in {\it every} characteristic) for degrees up to eight  and several other situations. 
 
\end{abstract}

\maketitle
\section{Introduction}

Fix an algebraically closed field $k$ of prime characteristic $p>0$.

 Let $f$ be a homogeneous polynomial of degree $d$ over $k$.
  The $F$-pure threshold  is an  invariant of the degree $d$ hypersurface  in projective space cut out by $f$. It can be understood as a  measure of how far the hypersurface  is from  being Frobenius split, analogously to how the log canonical threshold  measures the failure to be log Calabi-Yau  (see \cite{SchwedeSmith}).

   For {homogeneous} polynomials in $n$ variables of degree $d$, it is known that the F-pure threshold satisfies
 \begin{equation}\label{bounds}
\fpt{f} \leq \min(\frac{n}{d}, 1).
 \end{equation}
 However, for specific values of $n, d $ (and $p$), this bound is not sharp. In this paper, we establish  a sharp upper bound for every value of $n$ and $d$---indeed, 
 we show the maximum value of the F-pure threshold is either $\lambda = \min(\frac{n}{d}, 1)$ or a {\it specific} truncation of  a 
 base $p$ expansion of $\lambda$; see Theorems \ref{generic2} and \ref{generic}.

 Our results can be framed in terms of the {\it F-pure threshold function}  on the parameter space $\mathbb P({\text{Sym}}^d(k^n)^*)$ of all polynomials of degree $d$ in $n$ variables (up to scalar multiple). 
 This function  is constant on locally closed sets, and defines a finite stratification of  $\mathbb P({\text{Sym}}^d(k^n)^*)$, as we explain in Section 3.  Our Theorem \ref{generic2}  provides a formula for  the $F$-pure threshold of the polynomials in the {\it largest} stratum---a dense open subset of $\mathbb P({\text{Sym}}^d(k^n)^*)$---in terms of $p,d,$ and $n$.    In particular, our formula shows  the $F$-pure threshold of a generic degree $d$ polynomial over of field of characteristic $p$ is often {\it strictly smaller} than the log canonical threshold  of a generic degree $d$ polynomial over $\mathbb C$, 
 which is always precisely $ \min(\frac{n}{d}, 1).$ 
   
  Our work is inspired by work of  Hernandez, N\'u\~nez-Betancourt, Witt and Zhang \cite{HNWZ},
   who investigated the F-pure threshold of homogeneous polynomials cutting out a {\it smooth} projective hypersurface.
In two variables,  they prove the F-pure threshold of  a reduced\footnote{Note that a homogeneous polynomial in two variables over an algebraically closed field  defines a smooth subscheme of $\mathbb P^1$ if and only if it has no repeated factors, {\it{ i.e.}} is reduced.}  polynomial is either $ \frac{2}{d}$ or
{\it some} truncation of the non-terminating base $p$ expansion of $\frac{2}{d}$; they also  list some  necessary conditions on the place of the 
truncation and raise the issue as to what extent these might be sufficient. We  show that, while their conditions are {\it not} sufficient for all $d$ and $p$ (see Example \ref{NotSufficient}),   there
  {\it always exist} reduced polynomials of degree $d$ (provided $p$ does not divide $d$ and $d\geq 4$) whose F-pure threshold is {\it some} truncation of  $ \frac{2}{d}$; See Theorem \ref{order}.

   Ideally, we would like to understand, 
   for each $n, d$ and $p$, the {\it distinct} values of the F-pure threshold as we range over all reduced polynomials in $\mathbb P({\text{Sym}}^d(k^n)^*)$ (or over those defining a reduced, or a smooth,   projective hypersurface).
     We can easily describe equations for the locally closed set $X_{\lambda}$ of polynomials whose $F$-pure threshold is exactly $\lambda$ (see \S 3), but it may be difficult to  determine whether or not  a particular $X_{\lambda}$ is   {\it non-empty,} and if so, whether its intersection with the locus of reduced polynomials is non-empty. 
     We make progress towards this task by showing that when  $p$ does not divide $d$,  there is 
{\it always} a non-empty stratum of reduced polynomials of degree $d$ whose $F$-pure threshold is strictly less than $\frac{2}{d}$; see Theorem \ref{order2}. In Conjecture \ref{conj}, we speculate more precisely that there is always a {reduced} polynomial of degree $d$ whose $F$-pure threshold is 
{\it exactly }  the $e$-th truncation of $\frac{2}{d}$,  for  $e$ the order of $p$ in the multiplicative group 
$(\mathbb Z/d\mathbb Z)^\times$ (or in $(\mathbb Z/\frac{d}{2}\mathbb Z)^\times$ if $d$ is even).  We prove this conjecture when $e=1$ or $2$; see  
  Remark  \ref{e1} and Theorem \ref{p1}.  
  
  In Section 7, we determine the $F$-pure thresholds for all homogeneous polynomials in two variables in degrees up to eight, in each and  every  characteristic $p$. These computations provide further support for Conjecture \ref{conj}.

Lower bounds on the $F$-pure threshold of {\it reduced} polynomials are also of interest.\footnote{Without restricting to reduced polynomials, the minimal $F$-pure threshold is $\frac{1}{d}$, achieved by  $x^d$, regardless of the characteristic or number of variables.} 
 Proposition \ref{LB} provides a
 lower bound for the $F$-pure threshold of a  reduced degree $d$ polynomial as the first non-vanishing 
truncation of the base $p$ expansion of $\frac{2}{d}$. This raises the question: {\it what is the $F$-pure threshold for the {\bf {smallest}} stratum containing a reduced polynomial in  the stratification on $\PSymN$ by $F$-pure threshold  in terms of $n, d, $ and  $p$? } Of course, it is also of interest to ask the same about polynomials defining a {\it smooth} hypersurface. We make some specific speculations about the answers in the final section of the paper.

\noindent
 {\bf Acknowledgements.} Thanks to Luis N\'u\~nez-Betancourt for several helpful remarks on this work. 
 In addition, the first author is grateful to the University of Jyv\"askyl\"a for its hospitality while she undertook this project.

    
  \section{Preliminaries on F-pure threshold}

Fix  a  field $k$ of  positive characteristic  $p$. For an ideal $I$ in a ring of characteristic $p$, the notation $I^{[p^e]}$ denotes the ideal generated by the $p^e$-th powers of the elements (equivalently, any set of generators) of $I$.

\begin{definition}\label{FPureDef} The $F$-pure threshold of $f\in k[x_1, \dots, x_n]$ (at the  maximal ideal $\mathfrak{m}=(x_1, x_2, \dots,  x_n) $) is the real number
$$
\fpt(f) = \sup \left\{\frac{N}{p^e}  \,\,\, \bigg| \,\,\, f^N \notin \mathfrak{m}^{[p^e]} \right\} =  \inf \left\{\frac{N}{p^e}  \,\,\, \bigg| \,\,\, f^N \in \mathfrak{m}^{[p^e]} \right\},
$$
where  $\mathfrak{m}^{[p^e]} $ denotes the Frobenius power $\langle x_1^{p^e}, \dots, x_n^{p^e}\rangle$ of $\mathfrak m$.
\end{definition}

Originally, the F-pure threshold was defined somewhat differently  by Takagi and Watanabe  \cite{takagi+watanabe.F-pure_thresholds} as the "threshold" $c$ beyond which the pair $(S, f^{c})$ 
fails to be $F$-pure (see also  \cite{HaraYoshida}).  The reformulation above  has evolved through    \cite{mustata+takagi+watanabe.F-thresholds} and \cite{blickle+mustata+smith.discr_rat_FPTs}.  See also \cite{benito+faber+smith.measuring_sing_with_frob}.

\medskip The F-pure threshold should be viewed as a "characteristic $p$ analog" of the log canonical threshold, with smaller values of the F-pure threshold corresponding to "more singular" objects. Many of its  properties are familiar from corresponding statements about  the log canonical threshold.  Some properties are summarized  below in our
 setting: 

\medskip
\begin{prop} \label{basic}  Let $f$ be a homogeneous form of degree $d>0$ over a field $k$ of characteristic $p>0$.  Then 
\begin{enumerate}\label{easy}
\item  $\fpt(f) \leq  1$.
\item For a monomial $f = x_1^{a_1}x_2^{a_2}\cdots x_n^{a_n}$, $\fpt(f) = \min(\frac{1}{a_1}, \, \frac{1}{a_2}, \, \dots, \, \frac{1}{a_n})$.
\item $\fpt(f)\in \mathbb Q$. 
\item If $f$ is in $n$ variables, then  $\frac{1}{d}\leq \fpt (f) \leq  \min(1, \frac{n}{d})$.
\item 
For any $r\geq 1$, we have $\fpt(f^r)= \frac{\fpt(f) }{r}.
$
\item  $f^n\in \m^{[p^e]} $ if and only if $\fpt(f) \leq \frac{n}{p^e}$.
\item For any field extension $k'$ of $k$, the $F$-pure threshold is independent of whether we view $f$ as a polynomial over $k$ or over $k'$. 
\end{enumerate}
\end{prop}

\begin{proof} The first property is an immediate consequence of the definition, and the second follows easily as well. The third can be found in 
\cite{blickle+mustata+smith.discr_rat_FPTs}. The fourth and fifth are straightforward to check (see \cite[2.2]{extremal}) and the sixth follows by combining (5) with \cite[4.4]{extremal}.
For a proof of the final statement, see  \cite[3.4]{extremal}.
\end{proof}

In light of (7) in Proposition \ref{basic}, there is little loss of generality in studying $F$-pure threshold over an algebraically closed field. We therefore assume our ground field $k$ is algebraically closed throughout this paper.


      
  \subsection{Basics on Truncations}\label{TruncationBasics}  
      
       Fix a positive real number $\lambda \leq 1$, and consider a base $p$ expansion
       \begin{equation}\label{basep}
     \lambda =  \frac{\alpha_1}{p} + \frac{\alpha_2}{p^2} + \frac{\alpha_3}{p^3} +  \cdots   
     \end{equation}
      where the "digits" $\alpha_i$  are integers satisfying $0\leq \alpha_i \leq p-1$ for all $i\in \mathbb N$.
     If $\lambda$ is a rational number of the form $\frac{a}{p^e}$ with $a\in \mathbb Z$, then there are two such expansions (\ref{basep}), one terminating after the $e$-th place, and the other {\it non-terminating}, 
     with all digits equal to $p-1$ after the $e$-th place. Otherwise, the expression (\ref{basep}) for $\lambda$ 
     is unique, with the sequence of digits  $\{\alpha_i\}$ neither eventually zero nor eventually the constant sequence $\{p-1\}$.

     Thus every real number has a   {\it {unique  non-terminating}} base $p$ expansion,  in which  the sequence $\{\alpha_i\}$ of digits is not eventually zero. 
   We will use the notation $\langle \lambda\rangle_e$ for the {\it truncation at the $e$-th spot} of the {\bf non-terminating}  base $p$ expansion:
  \begin{equation}\label{TruncationNotation}
    \langle \lambda\rangle_e :=  \frac{\alpha_1}{p} + \frac{\alpha_2}{p^2} + \frac{\alpha_3}{p^3} +  \cdots   + \frac{\alpha_e}{p^e} = \frac{N_e}{p^e}.
    \end{equation}
We gather some facts about such truncations for future reference. The proofs are left as easy  exercises.

      \begin{lemma}\label{TrunBasic}  With notation as in (\ref{TruncationNotation}), 
      \begin{enumerate}
      \item Each truncation $\langle \lambda\rangle_e$ is {\it strictly} less than $\lambda$;
           \item  The sequence of truncations $\left\{\langle \lambda\rangle_e\right\}$ is a non-decreasing sequence converging to $\lambda$;
      \item The numerator $N_e$ of $\langle \lambda\rangle_e = \frac{N_e}{p^e}$ is the unique integer $M$ such that 
      $$
      0 < \lambda - \frac{M}{p^e} \leq \frac{1}{p^e};$$
      \item For every $e\geq 1$, $N_{e+1} = pN_e+\alpha_e$ where $0\leq \alpha_e\leq p-1$;
      \item For all positive integers $j\leq e$, the "round-down"  $\lfloor \frac{N_e}{p^{e-j}}\rfloor $ equals $N_j.$
      \end{enumerate}
          \end{lemma}
      
    \begin{lemma} \label{n/dTrunc}
    Let $\lambda=\frac{n}{d}$ be a rational number where $d, n \in \mathbb N$, and let $\frac{N_L}{p^L}$ be the truncation of the non-terminating base   $p$ expansion of $\lambda$.
Then  {\it provided $p^L\lambda\not\in \mathbb Z$}, when we divide the integer $p^Ln$ by $d$,  
 \begin{itemize}
 \item the 
 integer $N_L$ is the quotient, and 
 \item the integer $1\leq np^L-dN_L<d$ is the remainder.
 \end{itemize}
If $p^L\lambda \in \mathbb Z$, then the quotient of $np^L$ by $d$ is the integer $N_L+1$ (with no remainder).
    \end{lemma}

  
  \section{The stratification by F-pure threshold}
   
Fix an algebraically closed ground  field $k$ of characteristic $p$.

The F-pure threshold can be viewed as a function  assigning to each polynomial  $f\in k[x_1, \dots, x_n]$
 the rational number $\fpt(f)$. In this section, we prove some basic facts about this function that are likely  known to experts but hard to find in the literature. At the same time, we set up notation needed in the subsequent sections.

Let ${\text{Sym}}^d(k^n)^*$ denote the $k$ vector space of degree $d$ homogeneous polynomials in $n$ variables, and  $\mathbb P({\text{Sym}}^d(k^n)^*)$ the projective space of such polynomials up to non-zero scalar multiple. 
For a positive real number  $\lambda$,  define the following two subsets of  
$\mathbb P({\text{Sym}}^d(k^n)^*): $
$$
X_{\lambda} = \{f \,\, |\,\, \fpt(f) = \lambda \}  \,\,\,\,\,\,\,\,\,\, {\text{and}} \,\,\,\,\,\,\,\,\,\,
 X_{\leq \lambda} = \{f \,\, |\,\, \fpt(f)\leq  \lambda \}.$$
Note that $X_{\lambda} \subseteq 
 X_{\leq \lambda}$ and that $X_{\leq \lambda_1} \subseteq X_{\leq \lambda_2}$ whenever $\lambda_1 <\lambda_2.$

  \begin{prop}\label{Stratification} 
  The set $ X_{\leq \lambda}$   is a Zariski closed subset of $\mathbb P({\text{Sym}}^d(k^n)^*)$.
    \end{prop}
  
 Before proving Proposition \ref{Stratification}, we record the following obvious corollary:
 
 \begin{cor} \label{FPTExactLocallyClosed}
  The set  $X_{\lambda}$  is the locally closed set    of $\mathbb P({\text{Sym}}^d(k^n)^*)$
$$
  X_{\lambda} = X_{\leq \lambda}\,  \bigcap \, \left(\bigcup_{e\in \mathbb N} \mathbb P({\text{Sym}}^d(k^n)^*)\setminus X_{\leq \langle \lambda\rangle_{e}}\right)
$$
  where $ \langle \lambda\rangle_{e}$ is the truncation of the non-terminating base $p$ expansion of $\lambda$ at the $e$-th spot.
   \end{cor}
 
 \begin{proof}
 This follows immediately from Proposition \ref{Stratification} and the fact that the sequence $ \langle \lambda\rangle_{e}$ converges to $\lambda$  
 with every truncation $ \langle \lambda\rangle_{e}$ {\it strictly} less than $\lambda$. 
 \end{proof}

\begin{notation}\label{NotationForSym}
 Let $\mathcal{I}_d$ be the set of $n$-tuples ${\bf i}:=(i_1, \ldots, i_n)$ with $i_1, \ldots, i_n$ nonnegative integers such that  $i_1+\cdots + i_n=d$.
 The set of all degree $d$ monomials 
 $$
 \{x_1^{i_1}\cdots x_n^{i_n} \,\,\, | \,\,\,\, (i_1, \ldots, i_n)\in \mathcal{I}_d\}
 $$ is a basis for ${\text{Sym}}^d(k^n)^*$, the $k$-vector space of  all degree $d$ homogeneous polynomials in  $\{x_1, \dots, x_n\}$. We can thus
 write each $f\in {\text{Sym}}^d(k^n)^*$ uniquely as 
$$
f=\sum_{{\bf i} \in \mathcal{I}_d}  a_{{\bf i}}\, x_1^{i_1}\cdots x_n^{i_n},
$$
where $a_{\bf i}\in k$, and view the $a_{\bf i}$ as the coordinates for the space ${\text{Sym}}^d(k^n)^*$.   Thus the polynomial ring 
$$\mathcal{K}:=k[a_{\bf i}\, | \, {\bf i}\in \mathcal{I}_d]$$
 is the coordinate ring of ${\text{Sym}}^d(k^n)^*$, as well as the homogenous coordinate ring of the projective space $\mathbb P({\text{Sym}}^d(k^n)^*)$ of polynomials up to scalar multiple.  

   \medskip

  Consider  the {\it generic}  polynomial of degree $d$ in $n$ variables:
  \begin{equation}\label{GenPoly}
  f= \sum_{{\bf i} \in \mathcal{I}_d}  a_{{\bf i}}\, x_1^{i_1}\cdots x_n^{i_n} \in \mathcal K [x_1, \dots, x_n].
  \end{equation}
For each positive integer $N$, we have 
\begin{equation}\label{DefCC}
f^N=\sum_{{\bf j}=(j_1, \dots, j_n) \in \mathcal{I}_{Nd}} {\mathcal C}_{N {\bf j}} \, x_1^{j_1}\cdots x_n^{j_n}
\end{equation}
with ${\mathcal C}_{N {\bf j}}\in \mathcal{K}$ denoting the coefficient of $x_1^{j_1}\cdots x_n^{j_n}$ in the expansion of $f^N$.
Each ${\mathcal C}_{N {\bf j}}$ is a (possibly zero)
polynomial in the coordinates $a_{\bf i}$ of degree $N$.  
The following description of the polynomials ${\mathcal C}_{N\bf j}$, which will not be needed until \S 4,  is easy to verify by expanding out  $f^N$ and gathering together the $x_1^{j_1}\cdots x_n^{j_n}$ terms: 
\end{notation}

  \begin{formula}\label{FormC} Fix a positive integer $N$ and ${\bf j} = (j_1, \dots, j_n) \in \mathcal I_{Nd}$. 
For the polynomial  ${\mathcal C}_{N {\bf j}}\in \mathcal K$ as defined in (\ref{DefCC}), we have 
\begin{equation}\label{coeff1}
{\mathcal C}_{N {\bf j}}=\sum_{\bf k} \left(\frac{N!}{\prod_{{\bf i} \in \mathcal{I}_d} k_{{\bf i}}!}\right) \prod_{{\bf i} \in \mathcal{I}_d} a_{{\bf i}}^{k_{\bf i}}
\end{equation}
where the summation is over all choices of  $|\mathcal{I}_d|$-tuples ${\bf k} = (k_{\bf i}\, | \, {\bf i} \in \mathcal{I}_d)$ satisfying 
\begin{equation}\label{requirement1} \sum_{{\bf i} =(i_1, \dots, i_n)  \in \mathcal{I}_d} k_{\bf i}=N, \,\,\,\,\,\,\, {\text{and}}\end{equation}
\begin{equation}\label{requirement2} 
  \sum_{{\bf i}  =(i_1, \dots, i_n)\in \mathcal{I}_d}  k_{\bf i}i_1 =j_1, \ \ \ \ \ \  \sum_{{\bf i}  =(i_1, \dots, i_n)\in \mathcal{I}_d}  k_{\bf i}i_2 =j_2, \ \ \ \ \ \   \dots \ \ \ \ \ \ \ \ \ \ \  \sum_{{\bf i}  =(i_1, \dots, i_n)\in \mathcal{I}_d}  k_{\bf i}i_n=j_n.\end{equation}
\end{formula}
The equations (\ref{requirement2}) can be viewed as an expression of homogeneity of  $\mathcal C_{N{\bf j}}$ with respect to certain non-standard gradings on the polynomial ring  $\mathcal K$. For example, the first equality in (\ref{requirement2}) says that $\mathcal C_{N\bf j}$ is homogeneous of degree $j_1$ with respect to  the grading where
 $a_{\bf i}$ (for ${\bf i } =(i_1, \dots, i_n)\in \mathcal I_d$) has degree $i_1$; this grading tracks the degree of $x_1$ in the expression (\ref{DefCC}).


  \begin{proof}[Proof of Proposition \ref{Stratification}]
  Fix $\lambda\in \mathbb R_{>0}$. Note that the sequence
  \begin{equation}\label{decreasingSeq}
  \frac{\lceil p^e\lambda\rceil }{p^e} \,\,\,\,\,\,\, {\text{ for }} \,\,\, e \in \mathbb N
  \end{equation}
  is a non-increasing sequence of rational numbers  converging  to $\lambda$ from above.
  It follows immediately that 
 \begin{equation}\label{DecreasingIntersection}
 X_{\leq \frac{\lceil p^1\lambda\rceil }{p^1} }\supseteq  X_{\leq \frac{\lceil p^2\lambda\rceil }{p^2} } \supseteq X_{\leq \frac{\lceil p^3\lambda\rceil }{p^3} } \supseteq \ \cdots  \ \supseteq
  \bigcap_{e\in \mathbb N} X_{\leq \frac{\lceil p^e\lambda\rceil }{p^e} } =  X_{\leq \lambda}.
\end{equation}
  So to show that $  X_{\leq \lambda}$ is closed, it suffices to show that this is so when  $\lambda$  has  the form $\frac{N}{p^E}$.
  
  Recall that  $\fpt(f) \leq \frac{N}{p^E}$ if and only if $f^N\in \m^{[p^E]}$ where $\m=(x_1, \dots, x_n)$ (Proposition \ref{basic} (6)). Thus we need only observe that the  inclusion $f^N\in \m^{[p^E]}$
is a closed condition on  $ \PSymN$.  Such statements are well-known, but  we  do so by finding explicit defining polynomials in $\mathcal K$ because we will need to refer to them later. 
  
 For a polynomial $ \sum_{{\bf i} \in \mathcal{I}_d}  a_{{\bf i}}\, x_1^{i_1}\cdots x_n^{i_n} \in {\text{Sym}}^d(k^n)^*$, we
consider when 
\begin{equation}\label{cond1}
  f^N = \sum_{{\bf j}\in \mathcal I_{Nd} }\mathcal C_{N \bf j} x^{\bf j} \in (x_1^{p^E}, \dots, x_n^{p^E}).
  \end{equation}
 Since $(x_1^{p^E}, \dots, x_n^{p^E})$ is a monomial ideal, the inclusion (\ref{cond1}) holds if and only if 
 the polynomial 
\begin{equation}\label{NotZero}
 \mathcal C_{N\bf j } x^{\bf j} \in (x_1^{p^E}, \dots, x_n^{p^E}) \,\,\,\,\,\,\,\,\forall \,\, {\bf j}\in \mathcal I_{Nd}.
 \end{equation}
 This being automatic if $x^{\bf j}\in \m^{[p^E]}$, we see that  (\ref{cond1}) holds for a polynomial with coefficients $[a_{\bf i}] \in \PSymN$ if and only if 
 \begin{equation}\label{cond2}
 \mathcal C_{N \bf j }(a_{\bf i}) = 0 \,\,\,\,\,\,\, \forall \,\,\, {\bf j} =(j_1, \dots, j_n) \in \mathcal I_{Nd} \,\,\,\, {\text{such that }} \,\,\,\,j_k \leq p^E-1\,\,\, {\text{for}} \,\,\, k=1, 2, \dots n.
  \end{equation}
  Thus $X_{\leq \frac{N}{p^E}}$ is the closed set of $\PSymN$ defined by the polynomials (\ref{cond2}).
    \end{proof}
    
    For future reference, we highlight the final sentence of the proof of Theorem \ref{Stratification}:
    \begin{cor}\label{cond3}
 With notation as above, the locus in ${\text{Sym}}^d(k^n)^*$ of polynomials $f$ with  $\fpt(f) \leq \frac{N}{p^E}$ is the Zariski closed set 
 \begin{equation}\label{cond4}
  \mathbb V(\{\mathcal C_{N \bf j} \,\,\,\, | \,\,\, {\bf j} = (j_1, \dots, j_n) \in \mathcal I_{Nd} 
\,\,\,with \,\,\,  j_k < p^E\,\,\,  \forall k\}) \subset {\text{Sym}}^d(k^n)^*.
\end{equation}
In particular,   $X_{\leq \frac{N}{p^E}} = \mathbb P({\text{Sym}}^d(k^n)^*)$ if and only if  all
$\mathcal C_{N \bf j}$ whose indices ${\bf j}=(j_1, \dots, j_n)$ satisfy ${\bf j}<p^E$ component-wise  are {\it  zero} as polynomials in $\mathcal K$.
\end{cor}

  \medskip
 In light of Sato's  proof  of the ACC conjecture for F-pure threshold (see \cite[4.4]{BMS09}) we see that the 
 parameter space $\mathbb P({\text{Sym}}^d(k^n)^*)$  of degree $d$ polynomials in $n$ variables over a fixed field $k$ of characteristic $p$ admits a {\it finite} stratification by F-pure threshold: 
 \begin{cor} \label{FiniteStrat} 
 There is a finite descending chain of closed sets of  $\mathbb P({\text{Sym}}^d(k^n)^*)$,
\begin{equation}\label{strat2}
 X_{\leq \lambda_1} \subsetneq X_{\leq \lambda_2}  \subsetneq \cdots \subsetneq X_{\leq \lambda_T} = \mathbb P({\text{Sym}}^d(k^n)^*),
 \end{equation}
 where each  $X_{\leq \lambda_i} \setminus  X_{\leq \lambda_{i-1}} $ consists of polynomials of F-pure-threhold {\it exactly} $\lambda_i.$
 \end{cor}

 \begin{proof}
 The point is that there are only {\it finitely many}   {distinct} values of F-pure thresholds for a polynomial of degree $d$ in $k[x_1, \dots, x_n]$, in which case the statement follows immediately from Proposition \ref{Stratification}. 
 The finiteness holds because  the set of F-pure-thresholds for polynomials in  $\mathbb P({\text{Sym}}^d(k^n)^*)$ satisfies both the {\it ascending} and the {\it descending} chain conditions (ACC and DCC):  there is no strictly ascending,  nor any strictly descending sequence,  of F-pure thresholds among polynomials (in $n$ variables) of degree $d$ over $k$. The ACC condition is \cite[Thm 1.2]{Sato}, whereas the DCC condition follows from the Noetherianness of $\mathbb P({\text{Sym}}^d(k^n)^*)$. \end{proof}

 Corollary \ref{FiniteStrat} implies that there is a {\it maximal} value for the $F$-pure threshold of a homogeneous polynomial of degree $d$ in $n$ variables, and that the set of polynomial achieving this maximal $F$-pure threshold is a dense open set of $\PSymN$.
 We call this maximal $F$-pure threshold the {\bf generic} $F$-pure threshold for $\PSymN$. There are always reduced polynomials (and even polynomials defining a smooth projective hypersurface) achieving the generic $F$-pure threshold, since 
 these loci are also  dense and open in $\PSymN$.


\section{The largest stratum}

Fix an algebraically closed base field $k$ of characteristic $p$.

A basic question is: {\it what is the F-pure-threshold of the generic polynomial of degree $d$ in $n$ variables over $k$?} That is, what is the maximal F-pure threshold among all  polynomials  in $ \mathbb P({\text{Sym}}^d(k^n)^*)$?

  In this section, we answer this question, showing that the
generic F-pure threshold will always be either $\min(1, \frac{n}{d})$ or a truncation of its non-terminating base $p$ expansion at some particular spot. 
When $n\geq d$, this is easy:

\begin{prop}
For $n\geq d$,  a generic polynomial in $\PSymN$ has F-pure threshold $1$.
\end{prop}
  
  \begin{proof}
  By Proposition \ref{basic}, $\fpt(f)\leq 1$ for all $f\in \PSymN$, whereas the  monomial $x_1x_2\cdots x_d$ has F-pure threshold $1$. 
  \end{proof}

  The case where $d> n$ is covered by the next result:
     
  \begin{theorem}\label{generic2} Fix a degree $d\geq n$. The value of the F-pure threshold of a generic element in $\PSymN$--- or equivalently, 
   the maximum value of the F-pure threshold  of any homogenous polynomial of degree $d$ in $n$ variables--- is the truncation $\langle \frac{n}{d}\rangle_L$ at  the $L$-th spot of $\frac{n}{d}$,
   where $L$ is the smallest positive number such that
   \begin{enumerate}
   \item[(a)]    $p^L\frac{n}{d}$ is not an integer; and
    \item[(b)]   the remainder when 
  $np^L$ is divided by $d$ is strictly less than $n$.  
  \end{enumerate}
  If no such $L$ exists, then the generic value of the F-pure threshold is 
$\frac{n}{d}$.
\end{theorem}

  \begin{remark}
  Note that  condition (a) in Theorem \ref{generic2} holds automatically {\bf unless},   writing $\frac{n}{d}$  in lowest terms as $\frac{a}{b}$,  we have that  $b=p^e$ with $e\leq L$. 
  \end{remark}

 \bigskip
For polynomials in two variables, Theorem \ref{generic2} simplifies:
   
  \begin{cor}\label{genericCOR} Fix a degree $d\geq 2$. 
 The maximum F-pure threshold of a polynomial in $\text{Sym}^d(k^2)^*$   is equal to the truncation $\langle\frac{2}{d}\rangle_e$  of the non-terminating base $p$ expansion of  $\frac{2}{d}$, where 
$e$ is the  {\it smallest} positive value of $e$ such that  $2p^e\equiv 1 \mod d$.  \\ If no such $e$ exists, then the maximum F-pure threshold is $\frac{2}{d}$.
  \end{cor}
  
 \begin{example} The F-pure threshold of a generic polynomial of even degree in two variables is $\frac{2}{d}$. Indeed, if $d$ is even,  then there is no $e$ such that $2p^e\equiv 1\mod d$. 
 \end{example}

  \begin{remark} More  generally, if $p\equiv 1 \mod d \ $ (and $d\leq n$), Theorem \ref{generic2} implies that the generic polynomial of degree $d$ in $n$ variables over a field of characteristic $p$  has $F$-pure threshold $\frac{n}{d}$. In particular, fixing $n$ and  $d$ and letting $p$ vary, there are infinitely many values of $p$ for  which the $F$-pure threshold of the generic polynomial of degree $d$ in $n$ variables is $\frac{n}{d}$. This is expected, since the generic log canonical threshold over $\mathbb C$ in this case is $\frac{n}{d};$ see  \cite[Conj 3.6]{mustata+takagi+watanabe.F-thresholds}. 
    \end{remark}

  To prove  Theorem \ref{generic2}, we  first restate it as follows:
  
  \begin{theorem}\label{generic} Fix a degree $d\geq n$. 
   For each natural number $L$, let
  $
  \langle\lambda\rangle_L =  \frac{N_L}{p^L}
  $ be the truncation of the {\it non-terminating} base $p$ expansion of $\lambda$ at the $L$-th spot.  
  Then
  \begin{enumerate}
  \item[(i)]
  If $dN_L\leq np^L-n$ for all $L$, then the maximum value of the F-pure threshold of any polynomial in $n$ variables of degree $d$ is $\frac{n}{d}$. 
  \item[(ii)]
 Otherwise, the maximal F-pure threshold is equal to the smallest value of $\langle\frac{n}{d}\rangle_L$ such that $dN_L\geq np^L - n+1.$
  \end{enumerate}
  \end{theorem}
 
  \begin{proof}[Proof of the Equivalence of Theorems  \ref{generic2} and \ref{generic}]
  If $p^L\frac{n}{d}$ is not an integer, condition (b) in Theorem \ref{generic2} is an alternative way to write condition (ii) from Theorem \ref{generic} (see  Lemma \ref{n/dTrunc}).
  On the other hand, if $p^L\frac{n}{d}$ is an integer, then  in lowest terms
  $\frac{n}{d}=\frac{a}{p^e}$, where $L\geq e$, so 
   $$
\frac{n}{d} = \frac{N_{L}}{p^L} + \frac{1}{p^{L}}.
  $$
    In particular, multiplying by $p^Ld$, we have 
     $
dN_L=p^Ln-d \leq  p^Ln-n
  $
  for all $L\geq e$. So (i) in Theorem \ref{generic} holds for all $L\geq e$ in any case. 
  \end{proof}
    

   \medskip
   We now turn to the proof of Theorem  \ref{generic}.  
   The heart is the following non-vanishing result for the coefficient polynomials $\mathcal C_{N \bf j}$ as defined in (\ref{DefCC}):    
   
    \begin{Keylemma}\label{Non-ZeroC}
     Fix a degree $d$ and number of variables $n$. Suppose $M$ and $e$ are natural numbers such that  $\frac{M}{p^e}<\frac{n}{d} \leq 1$.
   Let  
     $$
     M_{j} = \left\lfloor \frac{M}{p^{e-j}}\right\rfloor.
     $$
     Then the following are equivalent 
     \begin{enumerate}
     \item[(a)]  $dM_j \leq np^j - n $ for all positive integers $j\leq e$;
     \item[(b)] For some index  ${\bf j}=(j_1, \dots, j_n)\in \mathcal I_{dM}$ with  $j_t\leq p^e-1$ for each $t=1, \dots, n$, the polynomial 
     ${\mathcal C}_{M\bf j}$  is non-zero as a polynomial in the coordinate ring $\mathcal K$ of $\PSymN$  (see   Notation \ref{NotationForSym}).
     \item[(c)] There is a homogeneous polynomial $g$ of degree $d$ in $n$ variables such that $\fpt(g)> \frac{M}{p^e}$.
     \end{enumerate}
     \end{Keylemma}

\begin{remark}\label{ValidM}
One situation where we will apply Lemma \ref{Non-ZeroC} is the following. Let $\langle \frac{n}{d} \rangle_e = \frac{N_e}{p^e}$ for some positive $e$. Then taking $M$ to be $N_e$, it follows that $M_{j}=N_j$ for all $j\leq e$.
\end{remark}

We will prove the Key Lemma \ref{Non-ZeroC} by induction on $e$. Lemma \ref{Non-ZeroCC} will provide the base case. For both the base case and the inductive step, we make use of the following easy fact:

\begin{lemma}\label{elementary}
Given nonnegative integers $a_1, \ldots, a_n, b$ with $b \le \sum_{i=1}^n a_i$, there exist nonnegative integers $b_1, \ldots, b_n$ such that $b_i \le a_i$ for all $i$, and $\sum_{i=1}^n b_i =b$.
\end{lemma}
\begin{proof}[Proof of Lemma \ref{elementary}] Induce on $b$. If $b=0$, take $b_i=0$ for all $i$. If $b>0$, then some $a_i>0$;  without loss of generality, say $a_1>0$.
 Consider the collection of non-negative integers $b-1, a_1-1, a_2, \dots, a_n$. By our inductive assumption, there exist $b_1', \ldots, b_n'$  such that $b_1'\leq a_1-1$ and 
 $b_i' \le a_i$ for all $i\geq 2$, with $\sum_{i=1}^n b_i' =b-1$. The proof is completed by taking $b_1=b_1'+1$ and $b_i=b_i'$ for $i\geq 2$.
 \end{proof}

       \begin{lemma}\label{Non-ZeroCC} With notation as above, 
if  $p>N$, then the polynomial ${\mathcal C}_{N {\bf j}}\in  \mathcal K$ is nonzero for all ${\bf j} \in \mathcal I_{Nd}$. 
\end{lemma}

\begin{proof}[Proof of Lemma \ref{Non-ZeroCC}]
Since  $N<p$, 
 the multinomial coefficients $\frac{N!}{\prod_{{\bf i} \in \mathcal{I}_d} k_{{\bf i}}!}$ are all nonzero, so looking at Formula \ref{FormC},  the conclusion follows once we show that for each ${\bf j} \in \mathcal I_{Nd}$, there exist choices of $(k_{\bf i} \, | \, {\bf i} \in \mathcal I_d)$  satisfying the requirements (\ref{requirement1}) and (\ref{requirement2}). In other words, we need to show that every monomial $x_1^{j_1} \cdots x_n^{j_n}$ with ${\bf j} = (j_1, \ldots, j_n) \in \mathcal I_{Nd}$ appears  in the expansion (\ref{FormC}) of $f^N$, where $f$ is the generic polynomial $\sum_{{\bf i}\in \mathcal I_d} a_{\bf i} x_1^{i_1}\cdots x_n^{i_n}$ of degree  $d$ in $n$ variables.

For this, we induce on $N$. If $N=1$, the claim is obvious. Let $N>1$, and fix $(j_1, \ldots, j_n) \in \mathcal I_{Nd}$. From Lemma \ref{elementary}, we know that we can choose $(i_1, \ldots, i_n) \in \mathcal I_d$ such that $(i_1, \ldots, i_n) \le (j_1, \ldots, j_n)$ component-wise. Since $\sum_{t=1}^n(j_t-i_t) = d(N-1)$, the inductive assumption implies that 
 the monomial $x_1^{j_1-i_1}\cdots x_n^{j_n-i_n}$ appears with non-zero coefficient in the expansion (\ref{FormC}) of $f^{N-1}$. Multiplying by $f$, therefore, there is a monomial 
$x_1^{j_1}\cdots x_n^{j_n}=(x_1^{i_1}\cdots x_n^{i_n})(x_1^{j_1-i_1}\cdots x_n^{j_n-i_n})$ appearing  in the expansion (\ref{FormC}) of $f\cdot f^{N-1} = f^N$.  This completes the proof.



\end{proof}

     \begin{proof}[Proof of Key Lemma \ref{Non-ZeroC}]    
     First note that (b) and (c) are equivalent by Corollary \ref{cond3} and Proposition \ref{basic} (6).

     We next show that (b) implies (a). If  $dM_j \geq np^j - n +1$ for some $j\leq e$, then $f^{M_j}\in (x_1^{p^j}, \dots, x_n^{p^j})$ for {\it all } polynomials $f$ of degree $d$ in $n$ variables. Raising to the $p^{e-j}$-power, also 
     $f^{p^{e-j}M_j}\in (x_1^{p^e}, \dots, x_n^{p^e})$. So since $M\geq p^{e-j}M_j$, we have $f^{M}\in (x_1^{p^e}, \dots, x_n^{p^e})$ for all $f\in {\text{Sym}}^d(k^n)^*$.  But now Corollary \ref{cond3}  implies that the polynomials 
 ${\mathcal C}_{M\bf j}$ are zero for all  ${\bf j}=(j_1, \dots, j_n)\in \mathcal I_{dM}$ with all components $j_t\leq p^e-1$.
     
   Finally, we prove (a) implies (b)  by induction on $e$. First suppose that $e=1$. In this case, (a) implies $M_1 =M <p$, so Lemma \ref{Non-ZeroCC} ensures that {\it all} 
    ${\mathcal C}_{M \bf j}$ are non-zero.  
 
For the inductive step, assume that $e \ge 2$. Let $\overline{M} = \lfloor\frac{M}{p}\rfloor$, and  $\overline{M_j} = \lfloor\frac{\overline M}{p^{(e-1)-j}}\rfloor$. Observe that 
\begin{enumerate}
\item[(i)] $\frac{\overline{M}}{p^{e-1}} \leq \frac{M}{p^e} < \frac{n}{d} \leq 1$
\item[(ii)]  $d\overline{M_j} \leq d \lfloor \frac{M}{p^{e-j}} \rfloor= dM_j \leq np^j-n$.
\end{enumerate}
Therefore, we can apply the inductive hypothesis to the natural numbers $\overline M$, $e-1$ to conclude that there is some non-zero polynomial ${\mathcal C}_{\overline{M} \bf J}\in \mathcal K$ 
with   index ${\bf J } = (J_1, \dots, J_n)\in \mathcal I_{d\overline{M}}$ and all  components
$J_t \leq p^{e-1}-1$.  Writing $M=p\overline{M} + \alpha$ where $0\leq \alpha<p$,
with this choice of index $\bf J$, therefore, the  products
\begin{equation}\label{one-term}
({\mathcal C}_{\overline{M} \bf J})^p{\mathcal C}_{\alpha {\bf K}}
\end{equation}
are non-zero polynomial for {\it all} choices of ${\bf k} = (k_1, \dots, k_n)\in \mathcal I_{d\alpha}$. 

Now let us look at the polynomials ${\mathcal C}_{M \bf L}$.
Fix an index ${\bf L } \in \mathcal I_{dM}.$  Expanding  $f^{M} = (f^{\overline{M}})^pf^{\alpha}$, we  see that  
\begin{equation}\label{summation}
\mathcal C_{M {\bf L}} = \sum_{{\bf j}\in \mathcal I_{d\overline{M}}} (\mathcal C _{\overline{M} {\bf j}})^p \mathcal C_{\alpha \ {\bf L}-p{\bf j}}.
\end{equation}
Here, we use the convention that $\mathcal C_{\alpha{\bf L}-p{\bf j}}=0$ if ${\bf L}-p{\bf j}$ has a negative component.  Note that the non-zero polynomial
  $ ({\mathcal C}_{\overline{M} \bf J})^p{\mathcal C}_{\alpha\ {\bf L-pJ}}$ as in (\ref{one-term}) would be one of the terms in the sum (\ref{summation}) {\it provided that } ${\bf L} \ge p{\bf J}$ component-wise.

We now claim that there is a choice of ${\bf L}=(L_1, \ldots, L_n) \in \mathcal I_{dM}$ with every component $L_i\leq p^e-1$ and  ${\bf L} \ge p{\bf J}$ component-wise (with $\bf J\in \mathcal I_{d\overline{M}}$ the special index chosen above). Equivalently, we claim  there is a choice of ${\bf K}=(K_1, \ldots, K_n) \in \mathcal I_{d\alpha}$ with components $K_i \le p^e-pJ_i-1$ for each $i$.
To see this, note that the assumption that $dM \le np^e-n$ gives 
$$
d\alpha =dM-pd\overline{M} \le np^e-n - \sum_{i=1}^np J_i =\sum_{i=1}^n (p^e-1-pJ_i).
$$
The existence of the claimed index $\bf L$ now follows from Lemma \ref{elementary}.



We will complete the proof of the Key Lemma  by showing that for this choice of $\bf L$,  the polynomial ${\mathcal C}_{M\ \bf L} $ is non-zero.
Because the non-zero polynomial $({\mathcal C}_{\overline{M}\ \bf J})^p{\mathcal C}_{\alpha {\bf K}}$ is one of the terms in the summation  (\ref{summation}), 
it suffices to show it  does not cancel in the sum.

To this end, fix  a monomial $\prod_{{\bf i} \in \mathcal I} a_{\bf i}^{l_{\bf i}}$
that has nonzero coefficient in $\mathcal{C}_{\overline{M}  {\bf J}}$, and  a monomial  $\prod_{{\bf i} \in \mathcal I} a_{\bf i}^{k_i}$ that has nonzero coefficient in 
${\mathcal C}_{\alpha \bf K}$ (where ${\bf K}={\bf L}-p{\bf J}.$) 
Then $\prod_{{\bf i} \in \mathcal I} a_{\bf i}^{pl_i + k_i}$ has nonzero coefficient in $(\mathcal{C}_{\overline{M} {\bf J}})^p\mathcal{C}_{\alpha {\bf K}}$.

We claim that the monomial $\prod_{{\bf i} \in \mathcal I} a_{\bf i}^{pl_i + k_i}$ does not appear in any other term  
$(\mathcal{C}_{\overline{M} {\bf j}})^p\mathcal{C}_{\alpha\ {\bf L}-p{\bf j}}$ for ${\bf j}\ne {\bf J}$, and therefore it cannot get canceled in the summation (\ref{summation}). 
To check this claim, assume, on the contrary,   that 
$\prod_{{\bf i} \in \mathcal I} a_{\bf i}^{pl_{\bf i} + k_{{\bf i}}}=\prod_{{\bf i} \in \mathcal I} a_{\bf i}^{pl'_{\bf i} + k'_{\bf i}}$ where $\prod_{{\bf i} \in \mathcal I} a_{\bf i}^{l'_{\bf i}}$ is a monomial appearing in $\mathcal{C}_{\overline{M} {\bf j}}$ and $\prod_{{\bf i} \in \mathcal I} a_{\bf i}^{k'_{\bf i}}$ is a monomial appearing in $\mathcal C_{\alpha \, {\bf L}-p{\bf j}}$. Note that $0\le k_{\bf i}, k'_{\bf i}\le p-1$ for all ${\bf i} \in \mathcal I$ (since $\sum_{{\bf i} \in \mathcal I} k_{\bf i} =\sum_{{\bf i} \in \mathcal I}k'_{\bf i}=\alpha \le p-1$), so $l_{\bf i}p+k_{\bf i}=l_{\bf i}'p+k'_{\bf i}$ implies $l_{\bf i}=l'_{\bf i}$ and $k_{\bf i}=k_{\bf i}'$ for all ${\bf i} \in \mathcal I$. 
We note that the monomial $\prod_{{\bf i} \in \mathcal I} a_{\bf i}^{l_{\bf i}}$ can only occur as a term in $\mathcal{C}_{\overline{M} {\bf j}}$ if 
$$\sum_{{\bf i} \in \mathcal I} k_{\bf i} i_1 = j_1,\,\,\,\,  \ldots, \,\,\,\, \sum_{{\bf i} \in \mathcal I} k_{{\bf i}} i_n=j_n.$$ In other words, there is a unique ${\bf j} $ such that $\prod_{{\bf i} \in \mathcal I} a_{\bf i}^{l_i}$ occurs as a term in $\mathcal{C}_{\overline{M} {\bf j}}$, 
so we must have ${\bf j} = {\bf J}$. This completes the proof of the Key Lemma.
\end{proof}


 Before proving Theorem \ref{generic}, we need one more lemma\footnote{In two variables, the lemma follows from \cite[4.4]{HNWZ}.}:
 
 \begin{lemma}\label{TruncAlways}
 Fix $d\geq n$. If the  F-pure threshold of the generic homogeneous polynomial of degree $d$ in $n$ variables is {\it not} $\frac{n}{d}$, then it is some truncation of the non-terminating base $p$ expansion of $\frac{n}{d}$.
 \end{lemma}
 
 \begin{proof}[Proof of Lemma \ref{TruncAlways}]
 The F-pure threshold of any $f\in \PSymN$ cutting out a {\it smooth} hypersurface in $\mathbb P^{n-1}$ is either $\frac{n}{d}$ or a smaller rational number with denominator a power of $p$ by Theorem 3.5 of \cite{HNWZ}. In particular, we can write the F-pure threshold of the {\it generic} hypersurface  in the form
 $$
\frac{ N_e-E}{p^e}
$$ for some non-negative integer $E$ and some truncation $\frac{N_e}{p^e}= \langle \frac{n}{d}\rangle_e$ of the non-terminating  base $p$ expansion of $\frac{n}{d}$.
We claim that $E$ must be $0$, which will prove the lemma.
 
 For this, note that since $\frac{N_e}{p^e}<\frac{n}{d}$,  we have  $dN_e<np^e$ and so if $E\geq 1$, 
 $$
 d(N_e-E) < np^e - dE\leq np^e - d \leq np^e-n.
 $$
  Let $M=N_e -E$ and note that
$$
 \lfloor \frac{M}{p^{e-j}}\rfloor= \lfloor \frac{N_e-E}{p^{e-j}}\rfloor \leq  \lfloor \frac{N_e}{p^{e-j}}\rfloor  = N_j \leq np^j-n
 $$ 
 for all positive $j\leq n$. So by the Key Lemma \ref{Non-ZeroC}, there is a some $g\in \PSymN$ such that $g^{N_e-E}\not\in \m^{[p^e]}$. Such  $g$ has F-pure threshold strictly larger than
 $\frac{N_e-E}{p^e}$ by Proposition \ref{basic} (6). In other words, $\frac{N_e-E}{p^e}$ can not be the maximal  (generic) F-pure threshold. This contradiction forces $E=0$.
So the generic  F-pure threshold is $\frac{N_e}{p^e}$, a truncation of $\frac{n}{d}$.
  \end{proof}

 We  now  have all the pieces in place to prove the theorem.
   
\begin{proof}[Proof of Theorem \ref{generic}]  
Suppose $dN_L\leq np^L-n$ for all $L$. Assume, on the contrary, that the maximal value of $\fpt(f)$ for $f\in \PSymN$ is strictly less than $\frac{n}{d}$.
Because there are only finitely many values of the F-pure threshold for polynomials in $ \PSymN$, and the truncations $\frac{N_e}{p^e}$ approach $\frac{n}{d}$ from below, 
 there must be some $e>0$ such that 
$\fpt(f) < \langle \frac{n}{d} \rangle_L$ for all $f\in \PSymN$ and  all $L\geq e$.
This says $f^{N_L}\in \m^{[p_L]}$ for {\it all} $f\in \PSymN$, and so the polynomials $\mathcal C_{N_L \bf j}$ are all zero for all large $L$ and all indices $\bf j$ component-wise less than $p^L$ (by Corollary \ref{cond3}). This contradicts 
the Key Lemma \ref{Non-ZeroC}.

Now suppose that $L$ is minimal such that $dN_L\geq dn-n+1$. This inequality forces $f^{N_L}\in m^{[p^L]} $ for all $f\in \PSymN$, so that the maximal F-pure threshold is bounded above by $\frac{N_L}{p^L}$, which is strictly less than $\frac{n}{d}$.  By Lemma \ref{TruncAlways}, the maximal F-pure threshold is therefore some truncation of $\frac{n}{d}$, so it is either the $L$-th truncation, or some earlier truncation.
By the minimality of $L$, the Key Lemma ensures that there is some $g\in \PSymN$ such that  $g^{N_{L-1}}\not\in \m^{[p^{L-1}]}$. So $\fpt(g)>\frac{N_{L-1}}{p^{L-1}}$ by Proposition \ref{basic} (6). This says the maximal F-pure threshold is strictly larger than the $L-1$-st truncation. So the maximal F-pure threshold must be $\frac{N_L}{p^L}= \langle \frac{n}{d} \rangle_L$, as claimed. 

\end{proof}


\section{Non-Empty strata for reduced polynomials in two variables}

We now take a closer look at the  stratification by F-pure threshold for {\it reduced} polynomials  in two variables.
 
 The stratification of  Proposition \ref{FiniteStrat} can be sharpened for polynomials in two variables  by incorporating results due to Hernandez, N\'u\~nez-Betancourt, Witt and Zhang---namely, 
 we see that {\it every} F-pure threshold of a polynomial defining a smooth\footnote{Since a zero dimensional scheme is smooth over an algebraically closed field 
  if and only if it is reduced, reduced polynomials in two variables are the same as polynomials in two variables defining a smooth hypersurface in $\mathbb P^1.$}
 projective hypersurface is either $\frac{2}{d}$ or some truncation:
 
  \begin{cor} \label{TwoVar} Fix a degree $d\geq 2$. There exists a natural number $e$ such that 
  the variety   $\ \mathcal U_{sm} \subset \mathbb P({\text{Sym}}^d(k^2)^*)$  of {\it reduced}  polynomials
  of degree $d$ in two variables  is a disjoint union
\begin{equation}\label{strat}
\mathcal U_{sm} = 
\mathcal U_{\langle\frac{2}{d}\rangle_{1}} \bigsqcup  \ \mathcal U_{\langle\frac{2}{d}\rangle_{2}}  \bigsqcup  \ \mathcal U_{\langle\frac{2}{d}\rangle_{3} } \  \cdots \  \ \bigsqcup \ \mathcal U_{\langle\frac{2}{d}\rangle_{e} } \bigsqcup \ \mathcal U_{\frac{2}{d}}
  \end{equation}
  where\footnote{Here, the notation $\langle \frac{2}{d}\rangle_L$ denotes  the truncation  of the non-terminating base $p$ expansion of  $\frac{2}{d}$ at the $L$-th spot; see Subsection \ref{TruncationBasics}.}
  $$ \mathcal U_{\lambda }= \{f\in  \PSym \,\, | \,\,f\, {\text{ {\bf reduced}, and  }}  \fpt(f) =  \lambda\}
  $$ is the locally closed locus of reduced polynomials whose F-pure threshold is precisely $\lambda$.
  \end{cor}
    
 We emphasize that some of the loci $ \mathcal U_{\langle\frac{2}{d}\rangle_{L} }$ or $ \mathcal U_{\frac{2}{d}}$ 
in the  stratification (\ref{strat}) can be empty.   Indeed, Hernandez, N\'u\~nez-Betancourt, Witt and Zhang  \cite[Thm 4.4]{HNWZ},  proved the following {\it necessary} conditions on $L$ such that 
$ \mathcal U_{\langle\frac{2}{d}\rangle_{L} }$ is non-empty: writing $\frac{2}{d}$ in lowest terms as $\frac{a}{b}$ and assuming $\langle\frac{2}{d}\rangle_{L} \neq 0$, 
 \begin{enumerate}
 \item[(I)] If $p$ does not divide $b$, then $L$ is bounded above by the order of $p$ in the multiplicative group $(\mathbb Z/b\mathbb Z)^\times$;
 \item[(II)] If $p>b$, then for every natural number $e'< L$,  the integer $a$ is less than the remainder when the integer $ap^{e'}$   is divided by $b$.
   \item[(III)]  The remainder when $ap^L$ is divided by $b$ is a positive integer less than or equal to $b-a.$ 
 \end{enumerate}
 A natural question\footnote{See the discussion in \cite{HNWZ} about "minimal lists" of conditions; eg Remark 4.8 and the subsequent examples.} is the {\it sufficiency} of the above three conditions: given values of $d$  and $L$ satisfying (I), (II) and (III), {\it does there exist a reduced homogeneous polynomial $f\in k[x, y]$ of degree $d$ such that $\fpt(f)$ is the truncation of the non-terminating base $p$ expansion of $\frac{2}{d}$ at the $L$-th term?} 
 While this turns out to be the case for small values of $d$ (see \S \ref{Smalld}), 
 the next example shows that these three conditions are not sufficient in general:
 
 \begin{example}\label{NotSufficient} Fix any $p>7$.
 Let $d=2p-3.$ Note that $L=2$ satisfies conditions (I), (II) and (III) above.
 The second truncation $\langle\frac{2}{d}\rangle_2$  is $\frac{1}{p}+\frac{1}{p^2}$.  
 However, $\frac{1}{p}+\frac{1}{p^2} $  is in the interval $(\frac{1}{p}, \frac{1}{p-1})$, and hence can not be the value of any F-pure threshold by \cite[4.3]{BMS09}.
 \end{example}

  Indeed, fixing arbitrary values of  $n, d$ and $p$,  it is not even obvious whether or not there is {\it any } non-empty
  $ \mathcal U_{\langle\frac{2}{d}\rangle_{L} }$. The next result settles this issue:

  \begin{theorem}\label{order2} Fix a characteristic $p$ and degree  $d\geq 4$, not divisible by $p$. 
  Then  there always  exist {\bf  reduced} polynomials $f\in k[x, y]$ of degree 
$d$ whose F-pure threshold is a truncation $ \langle \frac{2}{d}\rangle_L$ of the non-terminating base $p$ expansion of $\frac{2}{d}$ at the $L$-th spot, for some $L$.
 In particular, at least one 
$\mathcal U_{\langle\frac{2}{d}\rangle_{L} }$ in the decomposition (\ref{strat}) is non-empty.
\end{theorem}

When $d<4$, Theorem \ref{order2} fails:

\begin{remark}
If $d=2$, then every reduced polynomial has F-pure threshold $1 = \frac{2}{d}$:   up to linear change of coordinates, such a polynomial is $xy$, so this follows from  Proposition \ref{basic} (b). In this case,  the stratification (\ref{strat})  becomes $\mathcal U_{sm} =\mathcal U_{\frac{2}{d}}$ and  all 
$\mathcal U_{\langle\frac{2}{d}\rangle_L}$ are empty.
\end{remark}

\begin{remark}
If $d=3$, then every reduced polynomial can be put in the form $xy(x+y)$ after a linear change of coordinates, so all reduced polynomials have the same $F$-pure threshold.
In this case, the decomposition (\ref{strat}) has only one stratum
$$\mathcal U_{sm} = \begin{cases} \mathcal U_{\frac{2}{3}} & {\text {if}} \,\,p=3\,\, {\text{ or }} \,\,p\equiv 1 \mod 3\\
\mathcal U_{\langle\frac{2}{3}\rangle_{1}} = \mathcal U_{\frac{2p-1}{3p}}& {\text {if}} \,\,  p\equiv 2\mod 3,
\end{cases}
$$
where the F-pure thresholds are computed using
 Corollary \ref{genericCOR}. (See also \cite{HNWZ}).
\end{remark}

We will prove Theorem \ref{order2} by proving the following  stronger statement:

\begin{theorem}\label{order} 
Fix a characteristic $p$ and degree  $d\geq 4$, not divisible by $p$. Let $e$ be the smallest positive integer such that  $d$ divides $2p^e-2$. There there exists a reduced homogeneous polynomial in $k[x, y]$ whose F-pure threshold is $ \langle \frac{2}{d}\rangle_L$  for some $L\leq e$.
\end{theorem}

Indeed,  we conjecture that, under certain obvious necessary conditions, the locus $\mathcal U_{\langle\frac{2}{d}\rangle_{e} }$ itself  (with $e$ as in Theorem \ref{order}) is always non-empty; See Conjecture \ref{conj}.

\begin{remark}
If $p$ divides $d$, Theorem \ref{order2} can fail. For example, if $d=p^t$ or $d=2p^t$, then no value of $L$ satisfies condition (III) above.
In these cases, there is only one stratum in the filtration (\ref{strat})---every reduced polynomial has F-pure threshold equal to $\frac{2}{d}$. 
\end{remark}

Theorem \ref{order} is only superficially stronger than 
Theorem \ref{order2} because condition (I) already gives a bound on the truncation place. Indeed, there are several equivalent ways to think about the bound in condition (I):

\begin{prop}\label{Ways} For a prime $p$ and integer $d\geq 2$ not divisible by $p$,  the following are equivalent statements about a positive integer $e$:
\begin{enumerate}
\item  The integer $e$ is the smallest positive integer such that  $d$ divides $2p^e-2$;
\item Writing the fraction $\frac{2}{d}$ as $\frac{a}{b}$ in lowest terms, $e$ is the order of $p$ in the multiplicative group $(\mathbb Z/b\mathbb Z)^\times$;
\item The $e$-th truncation $\langle\frac{2}{d}\rangle_e$  is the fraction $\frac{2p^e-2}{dp^e}$, and no earlier truncation satisfies  $\langle\frac{2}{d}\rangle_L=\frac{2p^L-2}{dp^L}$;
\item The non-terminating base $p$ expansion of $\frac{2}{d}$ begins repeating with period $e$ at the $e$-th place.
\end{enumerate}
\end{prop}

\begin{proof} We mostly leave this as an exercise, pointing the reader to the basic facts established in Lemmas \ref{TrunBasic} and \ref{n/dTrunc}. For (4), 
 write $\frac{2}{d}=\frac{N_e}{p^e} +R$ where  $R=\frac{1}{p^e}\beta$ for some $\beta\leq 1$ (so $R$ is the remaining part of the base $p$ expansion).  Multiplying by $dp^e$, we have 
$2p^e=dN_e+dp^eR$, so solving for $R$ we see that $R=\frac{1}{p^e}\frac{2}{d}$.
\end{proof}


\bigskip

We now turn our attention to the proof of Theorem \ref{order}.
The idea is to show that for $e$  as in the theorem, the closed locus $X_{\leq \frac{N_e}{p^e}}$ is not entirely contained in the closed locus of non-reduced polynomials in $\PSym$.

\begin{notation}\label{SimplerNotation}
Simplifying  the notation defined earlier in (\ref{NotationForSym}), we let $a_0, \dots, a_d$ be the affine coordinates for ${\text{Sym}}^d(k^2)^*$ so that a polynomial of degree $d$ is written
$$
a_0x^d+a_1x^{d-1}y + \cdots+ a_dy^d.
$$
It makes sense to also denote the polynomials 
$\mathcal C_{N\, (j_1, j_2)}$, where $(j_1,j_2)\in \mathcal I_{dN}$ simply by 
$$
\mathcal C_{N j_2}, \,\,\,\,\,\,\, {\text{where }} j_2  {\text{  is an integer index satisfying  }} 0\leq j_2\leq dN, 
$$
as the value of $j_1$ is uniquely determined to be $dN-j_2$ in this case. So $\mathcal C_{N j}\in k[a_0,\dots, a_d]$ is the coefficient of $x^{dN-j}y^j$ in the expansion of $(a_0x^d+a_1x^{d-1}y + \cdots+ a_dy^d)^N$. Note that the degree $d$ is suppressed in our notation but the polynomials $\mathcal C_{Nj}$ very much depend on $d$.\end{notation} 

\noindent
{\bf The Discriminant.}
Recall that the {\bf discriminant}  $\Delta_d \in \mathcal K = k[a_0, \dots, a_d]$ is an irreducible homogeneous  polynomial defining the locus of non-reduced polynomials in  
${\text{Sym}}^d(k^2)^*$ (or $\PSym$). The discriminant $\Delta_d$ is invariant (up to scalar multiple) under the natural $GL_2(k)$ action on $\mathcal K$  induced by change of  coordinates $x$ and $y$. The reader can consult the first chapter of \cite{GKZ} for background on discriminants.
 
 We make use of the following  general lemma:

\begin{lemma}\label{Var}
Consider the natural $GL_2(k)$ action by change of coordinates on the affine space ${\text{Sym}}^d(k^2)^*$ of all polynomials of degree $d$ in two variables.
Let $\mathcal{C}\subseteq  k[a_0, \ldots, a_d]$ be a homogeneous radical  ideal stable under this  action, and let  $R$ be the quotient ring $k[a_0, a_1, \dots, a_d]/\mathcal C.$
If $\mathrm{dim}({R}/(a_{d-1}, a_d){R})=\mathrm{dim}({R})-2$, then there exist reduced polynomials $f =\sum_{i=0}^d a_ix^{d-i}y^i \in \mathbb V(\mathcal{C})$. Equivalently, the discriminant $\Delta_d \notin {\mathcal{C}}$.
\end{lemma}

\begin{proof}[Proof of Lemma \ref{Var}]  Let $\mathcal{X}\subseteq {\text{Sym}}^d(k^2)^*$ denote the  space of polynomials divisible by $x^2$, so that in the notation (\ref{SimplerNotation}), the set 
 $\mathcal X$ is the subvariety of ${\text{Sym}}^d(k^2)^*\cong k^{d+1}$ cut out by the equations 
$a_{d-1}=a_d=0$.
Note that $\mathcal X$ is a proper closed subset of the closed set $ \mathcal{N}=\mathbb V(\Delta_d)$  of all non-reduced polynomials of degree $d$.

Consider the regular map 
$$k \times \mathcal{X} \xrightarrow{\Phi} \mathcal{N}\,\,\,\,\,\,\,\,\, {\text{defined by}}\,\,\,\,\,\,\,\, \Phi(b, x^2g(x, y))=(x+by)^2g(x+by, y).$$
This is the action of the subgroup 
$$
G = \left\{
\begin{bmatrix} 1 & b \\ 
0 & 1 \end{bmatrix} \,\,\,\, | \,\, b\in k
\right\} \subset GL_2(k)
$$ on $\mathcal X$ by change of coordinates.  The image of $\Phi$ is the dense set of non-reduced polynomials divisible by $(x+by)^2$ for some $b \in k$---indeed, its complement in $\mathcal N$ is contained in the proper closed set  $\mathcal Y = \mathbb V(a_0, a_1)$ of polynomials divisible by $y^2$.  Thus the closure of $\mathrm{Im}({\Phi})$ is the full set $\mathcal N$ of non-reduced polynomials in ${\text{Sym}}^d(k^2)^*\cong k^{d+1}$.

Because the ideal $\mathcal C$ is stable under changes of coordinates, 
the map $\Phi$ restricts to a regular map $\tilde \Phi$ on the closed subsets given by intersection with  $\mathbb V(\mathcal C)\subset {\text{Sym}}^d(k^2)^*$:
$$
 k \times \tilde{\mathcal{X} }\xrightarrow{\tilde{\Phi}}  \tilde{\mathcal{N}}\,\,\,\,\,\,\,\,\, \,\,\,\,\,\,\,\,\, (b, x^2g(x, y))\mapsto (x+by)^2g(x+by, y),$$
where 
 $$\tilde{\mathcal{N}}:=\mathcal{N} \cap \mathbb V(\mathcal{C}) \,\,\,\,\,\,\,\, {\text{and}}\,\,\,\,\,\,\,\,\,\,  \tilde{\mathcal{X}}:=\mathcal{X}\cap \mathbb V(\mathcal{C}). $$
 It follows that the image of $\mathrm{Im}(\tilde{\Phi}) $ is  $\mathrm{Im}({\Phi})\cap \mathbb V(\mathcal C)$, so it is 
 dense in 
  $\tilde{\mathcal{N}}$. 

We want to show that the inclusion $\tilde{\mathcal{N}}\subset \mathbb V(\mathcal{C})$ is strict. For this, it suffices  if $\dim \tilde{\mathcal{N}} < \dim\mathbb V(\mathcal{C}) = \dim R$.
Since the image of $\tilde \Phi$ is dense in $\tilde{\mathcal N}$, it suffices to show  $\mathrm{dim}(\mathrm{Im}(\tilde{\Phi}))<\mathrm{dim}({R})$.
But since 
 $\mathrm{dim}(\mathrm{Im}(\tilde{\Phi})) \le \mathrm{dim} (k \times \tilde{\mathcal{X}}) = 1 + \mathrm{dim}(\tilde{\mathcal{X}})$, it is enough to show 
 $\mathrm{dim}(\tilde{\mathcal{X}}) \le \mathrm{dim}({R})-2$. This follows from our hypotheses, since $\tilde{\mathcal{X}}$ has coordinate ring ${R}/(a_{d-1}, a_d){R}$.
\end{proof}

\begin{proof}[Proof of Theorem \ref{order}]  
Fix $e$  to be the
smallest positive integer such that $d$ divides $2p^e-2$,   and observe that 
 $2p^e-2=N_ed$, where $\frac{N_e}{p^e}$ is the $e$-th truncation of $\frac{2}{d}$ (see Lemma \ref{n/dTrunc}).

We want to show that the closed set $X_{\leq \frac{N_e}{p^e}}$ in $ \PSym$  is not entirely contained in the locus 
$\mathbb V(\Delta_d)$ of non-reduced polynomials in  $\PSym$.
Note that because  $N_ed= 2p^e-2$, the  closed set $X_{\leq \frac{N_e}{p^e}}$ is defined by just {\it one} polynomial, $\mathcal C_{N_e p^e-1}$ in $k[a_0, \dots, a_d]$ (see Corollary \ref{cond3}). So to prove Theorem \ref{order}, we must show that 
\begin{equation}\label{WTSDisc}
\Delta_d \not\in \sqrt{(\mathcal C_{N_e p^e-1})}.
\end{equation}

Taking  $\mathcal C$ to be the principal ideal  $\sqrt{(\mathcal C_{N_e p^e-1})}$, note that $\mathcal C$ is invariant under linear changes of coordinates in $k[x, y]$, since the condition that $f^{N_e}\in \mathfrak{m}^{[p^e]}$ is. So by   Lemma \ref{Var},  
it  suffices to show that   $a_{d-1}, a_d$ forms a regular sequence on the ring $k[a_0, a_1, \dots, a_d]/ (\mathcal C_{N_e p^e-1})$. For this, it suffices to show that 
$a_{d-1}, a_d, \mathcal C_{N_e p^e-1}$ is a regular sequence on the polynomial ring $k[a_0, a_1, \dots, a_d]$. So to prove Theorem \ref{order}, we can prove that
\begin{equation}\label{suff}
\mathcal C_{N_e p^e-1}\not\in (a_{d-1}, a_d).
\end{equation}
But (\ref{suff}) is equivalent to the claim that there is {\it some} (non-reduced) polynomial  of the form 
\begin{equation}\label{x2g}
x^2(a_0x^{d-2}+a_1x^{d-3}y + \cdots + a_{d-2}y^{d-2}) 
\end{equation} which has F-pure threshold greater than $\frac{N_e}{p^e}$.
This is easy if $d$ is even: the polynomial $(xy)^{d/2} $  has F-pure threshold  $\frac{2}{d}$ by Proposition \ref{basic}, and $\frac{2}{d}$ is larger than any truncation. So we might as well assume $d$ is odd.

 Note that, for any $N$, 
\begin{equation}\label{x2g2}
\begin{aligned}
(x^2g)^N &= x^{2N}(a_0x^{d-2}+a_1x^{d-3}y + \cdots + a_{d-2}y^{d-2})^{N} \\
&=    x^{2N}\sum_{k=0}^{(d-2)N} {\mathcal C'_{N k}} x^{(d-2)N-k} y^k \\
 &= \sum_{k=0}^{(d-2)N} \overline{\mathcal C_{N k}} x^{dN-k}y^k,
 \end{aligned}
\end{equation}
where the coefficient of the monomial $x^{dN-k}y^k$ can be interpreted either as the polynomial 
$\mathcal C'_{N k}\in k[a_0, \dots, a_{d-2}]$ (arising from expanding out a generic polynomial of degree $d-2$)  {\bf or} as the 
image  $\overline{\mathcal C_{N k}}$ of $\mathcal C_{N k}$  modulo $(a_{d-1}, a_d)$ (arising from expanding out a generic degree $d$ polynomial).
 In particular, by 
 Lemma \ref{Non-ZeroCC}, we see that when $N<p$,  the polynomials
 \begin{equation}\label{smallN}
 \mathcal C_{Nk} \not\in (a_{d-1}, a_d) \,\,\,\,\,\,\,\,\,\, {\text { for all $k$ in the range}} \,\,\,0\leq k\leq (d-2)N.
 \end{equation}

 Now, to prove Theorem \ref{order},  we must show that $\mathcal C_{N_e p^e-1} \not\in (a_{d-1}, a_d)$. 
 First suppose $e=1$.  The assumptions $dN_1=2p-2$ and  $d$ odd imply that $N_1$ is even. Consider the degree $d$ polynomial $f=(xy)^{\frac{d-1}{2}}(x+y)$. We see that 
 $$
 f^{N_1}=(xy)^{\frac{N_1(d-1)}{2}}(x+y)^{N_1}
 $$
 has $(xy)^{p-1}$ term with coefficient the binomial coefficient  $\binom{N_1}{\frac{1}{2}N_1}$, which is non-zero since $N_1<p$. This says that $f^{N_1}\not\in (x^p, y^p)$, and so the $F$-pure threshold of $f$ is strictly greater than $\langle\frac{2}{d}\rangle_{1}=\frac{N_1}{p}$ as needed.
 
   Now assume $e\geq 2$.
  Writing $N_e=pN_{e-1}+\alpha$ where $\alpha\leq p-1$,
we have \begin{equation}\label{sum222}
\mathcal{C}_{N_ep^e-1} = \sum_{i} (\mathcal{C}_{N_{e-1} i})^p \mathcal{C}_{\alpha p^e-1-pi}
\end{equation} 
where the sum is taken over all non-negative $i$  satisfying $i\leq dN_{e-1}$ and $0\leq p^e-1-pi\leq d\alpha.$
To show that this polynomial is not in $(a_{d-1}, a_d)$, we  argue as in the proof of Lemma \ref{Non-ZeroC}: 
it suffices to find {\it one} valid choice of indices $i_0$ and $j_0$ with   $pi_0+j_0= p^e-1$ and both 
\begin{equation}\label{Suffices}
\mathcal{C}_{N_{e-1\, i_0}} \,\,\,\,\, {\text{and}}\,\,\,\,\,\,  \mathcal{C}_{\alpha j_0}  \,\,\,\,\, {\text{not in }}\,\,\,\,\,\, (a_{d-1},a_d).
\end{equation} Indeed,  reducing modulo $(a_{d-1}, a_d)$,  we see 
$(\overline{\mathcal{C}_{N_{e-1}\, i_0}})^p \overline{\mathcal{C}_{\alpha j_0}} $ is not the zero polynomial, and it admits some 
  monomial in $a_0, \ldots, a_{d-2}$ with nonzero coefficient which does not cancel in the summation (\ref{sum222}).
  This implies that $\overline{\mathcal{C}_{N_{e} p^e-1}}\neq 0$ and so $ \mathcal{C}_{N_{e} p^e-1}\not\in (a_{d-1}, a_d)$.
  
  To find the needed indices  $i_0$ and $j_0$, note that by (\ref{smallN}) above, since $\alpha<p$, the only restriction needed on $j_0$ is that $j_0\leq (d-2)\alpha$. To find $i_0$, we  use the following lemma, whose proof we defer until after finishing the proof of Theorem \ref{order}: 

\begin{lemma}\label{claimj} With notation as above, fix $d\geq 5$ and a positive integer $j\leq e-1$.
Then there exists an index $i$ with $dN_j-p^j +2 \le i \le p^j -1$, such that the polynomial $\mathcal{C}_{N_j i}\notin (a_{d-1}, a_d)$.
\end{lemma}

Using Lemma \ref{claimj}, we  choose $i_0$ such that 
$$
dN_{e-1}-p^{e-1} +2 \le i_0 \le p^{e-1} -1 \,\,\,\,\,\,\,\, {\text{and}} \,\,\,\,\,\,  \mathcal C_{N_{e-1\, i_0}}\not\in (a_{d-1}, a_d).
$$
Then 
$(\mathcal C_{N_{e-1 \, i_0}})^p\mathcal C_{\alpha\ p^e-i_0-1}$ is one of the terms of the sum (\ref{sum222}) 
  and  $\mathcal C_{\alpha \, p^e-i_0-1} \not\in (a_{d-1},a_d)$ provided that 
 \begin{equation}\label{NeedRange}
 0 \leq p^e-1-pi_0\leq (d-2)\alpha.
 \end{equation}
 
 The first inequality in (\ref{NeedRange}) follows from the fact that $i_0\leq p^{e-1}-1$. To get the second inequality in (\ref{NeedRange}), first use the fact that $i_0 \ge dN_{e-1}-p^{e-1}+2$ and  $pN_{e-1}=N_e-\alpha$,
to see that 
$$
 p^e-1-pi_0\leq  p^e-1-p(dN_{e-1}-p^{e-1}+2) = 2p^e-dN_e+d\alpha-2p-1.
 $$
Then using that   $dN=2p^e-2$, 
we see that 
$$
 p^e-1-pi_0\leq  d\alpha-2p+1 \leq (d-2)\alpha,
 $$
 so (\ref{NeedRange}) holds.
This completes the  proof of Theorem \ref{order}, once we have proved  Lemma \ref{claimj}.


\begin{proof}[Proof of Lemma \ref{claimj}]
  Fix $e>1$, else the statement is vacuous. We induce on $j$.
  
We may assume that $2p^j\not\equiv 1 \mod d$ for all $j\leq d$. For if  $2p^j\equiv 1 \mod d$, then Corollary \ref{genericCOR} ensures that the generic polynomial has $F$-pure threshold $\langle \frac{2}{d}\rangle_j$, which is less than or equal to $\langle \frac{2}{d}\rangle_e$, as needed.
Therefore, also $dN_j\leq 2p^j-2$ for all $j\leq e$ (by the opening lines of the proof of Corollary  \ref{genericCOR}.)
Further,  since $e$ is the {\it smallest} natural number such that $dN_e=2p^e-2$, we see also that 
\begin{equation}\label{TightBound}
dN_j\leq 2p^j-3 \,\,\,\,\,\,\,\,\, {\text{for all }} \,\,\, j \leq e-1.
\end{equation}

{\bf Base case:}  Assume $j=1$.  Note that $N_1\leq p-1$.  So by (\ref{smallN}), we know that $\mathcal C_{N_1 i}\not\in (a_{d-1}, a_d)$ for all non-negative values of $i\leq (d-2)N_1$.
We need only make sure that some such value of $i$ also satisfies 
 \begin{equation}\label{rangeNEW}
 dN_1-p+2\le i \le p-1.
 \end{equation}

Since $d> 4$,   using  (\ref{TightBound}),  we have $2N_1 < \frac{dN_1}{2}\le \frac{2p-3}{2}.$  Since $2N_1$ is an integer, this  implies $2N_1 \le p-2$.
Adding $(d-2)N_1$ to both sides and rearranging, we have
$$dN_1-p+2\le (d-2)N_1.$$
 So taking $i=\min(p-1, (d-2)N_1)$ produces an $i$ in the desired range, proving the base case.

{\bf Inductive step:} For $j\geq 2$, write
 $N_j=pN_{j-1}+\alpha$, where  $0 \le \alpha\le p-1$. 
We have $f^{N_j}=(f^{N_{j-1}})^p f^{\alpha}$, so 
$$
\sum_{k=0}^{dN_j} \mathcal{C}_{N_j k}x^{N_jd-k}y^k=\left(\sum_{i=0}^{dN_{j-1}} \mathcal{C}_{N_{j-1} i}^px^{(dN_{j-1}-i)p}y^{pi}\right) \left( \sum_{j=0}^{d\alpha} \mathcal{C}_{\alpha j}x^{\alpha d -j}y^j\right)
$$
and therefore
\begin{equation}\label{sumc}
\mathcal{C}_{N_j k}=\sum_{i, j} (\mathcal{C}_{N_{j-1} i})^p \mathcal{C}_{\alpha, j}
\end{equation}
where the summation is over all the choices of $0  \le i \le dN_{j-1}$ and $0 \le j \le d\alpha$ with $pi+j=k$.

If $I=pi_0+j_0$ for some $i_0$ and $j_0$ with $\mathcal{C}_{N_{j-1} i_0}, \mathcal{C}_{\alpha j_0} \notin (a_{d-1}, a_d)$, then 
$(\mathcal{C}_{N_{j-1} i_0})^p \mathcal{C}_{\alpha, j_0} \notin (a_{d-1}, a_d)$, because the ideal $(a_{d-1}, a_d)$ is prime. We claim that in this case, the polynomial
 $\mathcal{C}_{N_j I} $ can not be in the ideal $ (a_{d-1}, a_d)$ either. Indeed, working modulo $(a_{d-1}, a_d)$, we can argue as in the last part of the proof of Lemma \ref{Non-ZeroC}:  picking a monomial $a_0^{k_0}\cdots a_{d-2}^{k_{d-2}}$ that has a nonzero coefficient in $\mathcal{C'}_{N_{j-1} i_0}$ and a monomial  $a_0^{l_0}\cdots a_{d-2}^{l_{d-2}}$ that  has nonzero coefficient in  $\mathcal{C'}_{\alpha j_0}$, the monomial  $a_0^{pk_0+l_0}\cdots a_{d-2}^{pk_{d-2}+l_{d-2}}$ has nonzero coefficient in $(\mathcal{C}_{N_{j-1} i_0})^p \mathcal{C}_{\alpha j_0}$ and cannot cancel in the summation (\ref{sumc}).

From the inductive hypothesis, there is an $i_0$ in the range 
 \begin{equation}\label{picki2}
dN_{j-1}-p^{j-1}+2 \le i_0 \le p^{j-1}-1
\end{equation}
 such that $\mathcal{C}_{N_{j-1} i_0} \notin (a_{d-1}, a_d)$.
We want to find an integer $j_0$ such that $j_0+pi_0$ is in the interval $[dN_j-p^j+2, p^j-1]$ and the polynomial 
$\mathcal C_{\alpha\, j_0}\not\in (a_{d-1}, a_d)$. This last condition is equivalent to
 $j_0\leq (d-2)\alpha$ by (\ref{smallN}). 
 So we are looking for  an integer $j_0$ in the intersection of the intervals 
 \begin{equation}\label{pickj2}
 [0, (d-2)\alpha]\,\,\,\,\,\,\,\,\,\, {\text{and}} \,\,\,\,\,\,\,\,\,\,\,  [dN_j-p^j+2-pi_0, \, p^j-pi_0-1].
\end{equation}
The second interval is non-empty because $dN_j \le 2p^j -3$. 
Thus to ensure that the intervals  (\ref{pickj2})  have nonempty intersection, it suffices if
\begin{equation}\label{in1}
p^j-pi_0 -1 \ge 0, \,\,\,\,\,\,{\text{and}}
\end{equation}
\begin{equation}\label{in2}
dN_j-p^j+2-pi_0 \le (d-2)\alpha.
\end{equation}
Multiplying (\ref{picki2}) by $p$,  we have
\begin{equation}\label{in3}
pdN_{j-1} -p^j + 2p \le pi_0 \le p^j-p
\end{equation}
Inequality (\ref{in1}) follows from (\ref{in3}) since $pi_0 \le p^j-p < p^j-1$. To get inequality (\ref{in2}), 
observe
$$
\begin{aligned}
 pi_0 & \ge pdN_{j-1}-p^j+ 2p &{\text{from (\ref{in3})}}\,\,\,\,\,\,\,\,\,\, \,\,\,\,\,\,\,\,\,\,\,\,\,\,\,\,\,\,\,\,\,\,\,\,\,\, \,\,\,\,\, \\ 
 & \,\,\,\,\,\,\,\,\,\,=dN_j-d\alpha -p^j+2p  &{\text{because}}\,\,\,\, pN_{j-1}=N_j-\alpha \\ 
 & \,\,\,\,\,\,\,\,\,\,\,\,\,\,\,\,\,\,\,\,\ge dN_j-d\alpha -p^j+2\alpha +2  &\,\,\,\,\,\,\, {\text{because}}\,\,\,\, \alpha \le p-1\,\,\,\,\,\,\,\, \,\,\,\,\,\,\,\, \\ 
 &\,\,\,\,\,\,\,\,\,\,\,\,\,\,\,\,\,\,\,\,\,\,\,\,\,\,\,\,\,\, =dN_j-(d-2)\alpha -p^j+2.
\end{aligned}
$$
This concludes the proof of Lemma \ref{claimj}, and hence of Theorem \ref{order}.
\end{proof}
\end{proof}

\section{The existence of reduced polynomials  with specific F-pure thresholds}

A natural question is whether there is any specific truncation of $\frac{2}{d}$  that is {\it always} the F-pure threshold of a reduced polynomial in two variables.

For example, we speculate:
\begin{conjecture}\label{conj} Fix any characteristic $p$ and degree $d\geq 4$ not divisible by $p$. Let $e$ be the smallest positive integer\footnote{See Proposition \ref{Ways} for  some equivalent ways to describe the natural number $e$.} such that $2p^e\equiv2\mod d$. 
Then there is always a reduced polynomial in $k[x,y]$ 
whose F-pure threshold is equal to the truncation $$\langle\frac{2}{d}\rangle_e \,\,\,\,\,\,\, {\text{or equivalently, }}\,\, \frac{2p^e-2}{dp^e},$$ provided that $2p^L\not\equiv 1\mod d$ for any $L<e$.
\end{conjecture}

\begin{remark}
The condition that $2p^L\not\equiv 1\mod d$ for any $L<e$ is clearly {\it necessary}. Indeed,  if $2p^L\equiv 1\mod d$ for some value of $L$, then  Corollary  \ref{genericCOR} guarantees there is no polynomial in  $k[x, y]$  with
$F$-pure threshold larger than $\langle\frac{2}{d}\rangle_L$.
\end{remark}

\begin{remark}\label{e1}
Another way to think about Conjecture \ref{conj}:
there is a reduced polynomial in $k[x,y]$  whose $F$-pure threshold is $\langle \frac{2}{d}\rangle_e $, where $e$ is smallest natural number  such that $2p^e\equiv 1 {\text{  or  }}  2 \, \mod d$.
\end{remark}

Theorem \ref{order} implies immediately that Conjecture \ref{conj} holds when the value of $e$ is one.
Our next goal is to prove Conjecture \ref{conj} holds when  $e$ is two: 

\begin{theorem}\label{p1} Fix $p$ and $d\geq 4$ not divisible by $p$. Then if 
$2p^2 \equiv 2 \mod d$ and $2p\not\equiv 1\mod d$, then 
there exist {\bf reduced } $f \in k[x, y]$ of degree $d$ such that $\mathrm{fpt}(f)=\langle\frac{2}{d}\rangle_2=\frac{2p^2-2}{dp^2}$.
\end{theorem}

The key point in the proof of Theorem \ref{p1} is the following:

\begin{theorem}\label{codim2}
Fix a prime $p$ and degree $d\geq 4$ not divisible by $p$. Assume that the remainder, when $2p$ is divided by $d$ is at least three---that is, we can write $2p=Nd+r$ with $N$  non-negative and $3\leq r \leq d-1$.

Then the closed set $X_{\leq \frac{N}{p}} $ of polynomials in two variables of degree $d$ whose F-pure threshold is bounded above by $\frac{N}{p}$ has codimension at least two in $\PSym$.
\end{theorem}

\begin{proof} [Proof that Theorem \ref{codim2} implies Theorem \ref{p1}]
To prove Theorem \ref{p1},   it suffices to show that
$$
X_{\leq \frac{N_2}{p^2}}\not\subseteq \, X_{\leq \frac{N_1}{p}} \,\, \bigcup \,\,  \mathbb V(\Delta_d),
$$
where 
$\Delta_d$ is the degree $d$ discriminant (defining the locus of non-reduced polynomials in $\PSym$. Note that because $dN_2=2p^2-2$, 
the closed set $X_{\leq \frac{N_2}{p^2}}$  is defined by just one polynomial, $\mathcal C_{N_2\, p^2-1}$, so has dimension $d-1$.
So if 
$$
X_{\leq \frac{N_2}{p^2}}\subseteq \, X_{\leq \frac{N_1}{p}}\,\, \bigcup \,\, \mathbb V(\Delta_d),
$$ then assuming Theorem \ref{codim2},  no component of  the dimension $d-1$ subvariety $X_{\leq \frac{N_2}{p^2}}$ can be contained 
$X_{\leq \frac{N_1}{p}},$ because $X_{\leq \frac{N_1}{p}}$ has strictly smaller dimension. That is,  {\it every} component of $X_{\leq \frac{N_2}{p^2}} $ is contained in  $\mathbb V(\Delta_d),$ and so  $X_{\leq \frac{N_2}{p^2}}\subseteq  \mathbb V(\Delta_d)$. But  this  means that the $F$-pure threshold of every reduced polynomial is strictly greater than $\frac{N_2}{p^2}= \langle \frac{2}{d}\rangle_2$, contrary to  Theorem \ref{order}.
Thus Theorem \ref{p1} follows from Theorem \ref{codim2}.
\end{proof}

To prove Theorem \ref{codim2}, we will use the following technical lemma:

\begin{lemma}\label{genL1}
Fix  a prime $p$ and degree  $d\geq  4$.   Write   $2p = dN +  r$ where $N$ is non-negative and $3\leq r\leq d-1$. 
Suppose there exist indices $i, j $ such that
\begin{itemize}
\item  $0 \leq i, j <\frac{p}{N}$  
\item $d-1\leq i + j \leq d-2$.
\end{itemize}
Then for any degree $d$ polynomial of the form  $f=x^iy^jg$ with $g$ not divisible by $x$ or $y$, we have 
we have $f^N\notin (x^p, y^p)$.
\end{lemma}

To prove Lemma \ref{genL1}, we need the following "folklore" fact:

\begin{lemma}\label{colonLEM}
Let $h_1, \dots, h_n$ and  $g_1, \dots, g_n$ be homogeneous regular sequences in the polynomial ring 
 $k[x_1, \dots,  x_n]$. Then assuming that $D:=\sum_{i=1}^n\deg g_i- \sum_{i=1}^n\deg h_i>0,$
$$
(g_1, \dots, g_n): (h_1, \dots, h_t) \subseteq (g_1, \dots, g_n) + (x_1,\dots,  x_n)^{D}.
$$
\end{lemma}

\begin{proof}
Set $d=\sum_{i=1}^n\deg g_i.$ We first claim that for all $N\in \mathbb N$, 
\begin{equation}\label{stepHilb}
(g_1, \dots, g_n): (x_1,\dots,  x_n)^{N} \subseteq (g_1, \dots, g_n) + (x_1,\dots,  x_n)^{d-n+1-N}.
\end{equation}
Indeed, when $N=1$, this follows from the fact that 
$$k[x_1, \dots, x_n]/(g_1, \dots, g_n) \,\,\,\, {\text{and}} \,\,\, \,
k[x_1, \dots, x_n]/(x_1^{\deg g_1}, \dots, x_n^{\deg g_n})$$ have the same Hilbert function, and hence the same socle degree, $d-n$. 
For $N\geq 2$, the inclusion (\ref{stepHilb}) follows by induction on $N$.

Next, set  $d'=\sum_{i=1}^n\deg h_i.$  Note that $(x_1,\dots,  x_n)^{d'-n+1}\subseteq (h_1, \dots, h_t)$,  since,  again,  the socle degree of the quotient  $k[x_1, \dots, x_n]/(h_1, \dots, h_t)$ is $d'-n$. So 
$$
(g_1, \dots, g_n): (h_1, \dots, h_t) \subseteq 
(g_1, \dots, g_n) : (x_1,\dots,  x_n)^{d'-n+1}  \subseteq  (x_1,\dots,  x_n)^{d-d'},
$$
with the second inclusion obtained by applying  (\ref{stepHilb}) with $N=d'-n+1$.
\end{proof}

\begin{proof}[Proof of Lemma \ref{genL1}] 
Note that our assumptions imply that $g$ is either degree one or two.

Assume $f^N \in(x^p, y^p)$, or, equivalently, $g^N \in (x^{p-iN}, y^{p-jN})$. Let $N'$ denote the smallest positive integer satisfying $g^{N'} \in (x^{p-iN}, y^{p-jN})$. The assumption that $i, j <\frac{p}{N}$ guarantees that $(x^{p-iN}, y^{p-jN})$ is not the unit ideal.

\noindent
{\bf Case 1: $g$ reduced}. If  $g$ is reduced and not divisible by $x$ or $y$, the reader will easily verify that $x\frac{\partial g}{\partial x}, y\frac{\partial g}{\partial y}$ form a regular sequence consisting of forms of degree equal to $\deg  g$, which is 
$d-(i+j)$.

Applying the operators\footnote{A similar trick is used in \cite{HNWZ} as well as in \cite{BhattSingh}.} $x\frac{\partial }{\partial x}, y \frac{\partial}{\partial y}$ to the equation $g^{N'} \in (x^{p-iN}, y^{p-jN})$ yields 
$$
g^{N'-1} \in (x^{p-iN}, y^{p-jN}): (x\frac{\partial g}{\partial x}, y \frac{\partial g}{\partial y}),
$$
which  is contained in $(x^{p-iN}, y^{p-jN}) + \m^{2p-(i+j)N-2d+2(i+j)}$ by Lemma \ref{colonLEM}.
Since $g^{N'-1}\notin (x^{p-iN}, p^{p-jN})$, this implies
$$
\mathrm{deg}(g^{N-1}) \geq \mathrm{deg}(g^{N'-1}) \ge 2p -(i+j)N-2d+2(i+j).
$$
Thus, 
$
(N-1)(d-(i+j)) \ge 2p -(i+j)N-2d+2(i+j),
$
which simplifies to 
$$
Nd \ge 2p-d+(i+j).
$$
Recalling that $Nd=2p-r$, we obtain $i+j \le d-r\leq d-3$, which contradicts our hypothesis.

\noindent
{\bf Case 2: $g$ not reduced}. In this case, $g$ has degree two, and changing coordinates, we can assume that $f=x^iy^j(x+y)^2.$
Assume, on the contrary, that $f^N \in(x^p, y^p)$, then $(x+y)^{2N} \in (x^{p-iN}, y^{p-jN})$. Because  all binomial coefficients $\binom{2N}{\ell}$ are nonzero (note $2N<p$), this implies that
  every monomial $x^{\alpha}y^{\beta}$ with $\alpha+\beta=2N$ is in $(x^{p-iN},y^{ p-jN})$. Taking $(\alpha, \beta)=(p-iN-1, 2N-(p-iN-1))$, it follows that $2N-p+iN+1\ge p-jN$. 
 Rearranging, we get $(2+i+j)N\ge 2p-1$, and hence $dN=dN+r-1$, contradicting our assumption on $r$. 
\end{proof}

We are now ready to prove Theorem \ref{codim2}, and hence complete the proof of Theorem \ref{p1}.

\begin{proof}[Proof of Theorem \ref{codim2}]
We know that the closed set $X_{\leq \frac{N}{p}}$ is defined by the 
 polynomials 
\begin{equation}\label{Polys}
\mathcal C_{N\,p-1}, \mathcal C_{N\,p-2}, \dots, \mathcal C_{N\,p-(r-1)}
\end{equation}
in $k[a_0, a_1, \dots, a_d]$, using the notation of Subsection  \ref{SimplerNotation}.
If the ideal generated by the polynomials (\ref{Polys}) has height two or more, we are done: $X_{\leq \frac{N}{p}}$ has codimension two in $\PSym$. So assume the ideal has height one, in which case 
the polynomials (\ref{Polys}) have a common factor, since all height one primes in a polynomial ring are principal. Let $\mathcal C$ be the product of the irreducible common factors, and  observe that $\mathcal C$ has the following properties:
\begin{enumerate}
\item[(A)]\, $\mathcal C$ is homogeneous under the standard grading on $k[a_0, \dots, a_d]$;
\item[(B)] \, $\mathcal C$ is homogeneous under the  grading on $k[a_0, \dots, a_d]$ assigning degree $i$ to $a_i$;
\item[(C)]  \, $\mathcal C$ is invariant under interchanging $a_i$ and $a_{d-i}$ for all $i=0, 1, \dots d.$
\end{enumerate}
Indeed, (A) and (B) hold because each $\mathcal C_{N\, j}$ is homogeneous under the respective  grading (see the discussion immediately following (\ref{requirement2})), hence so are the irreducible divisors of each, as well as their product $\mathcal C$. And (C) holds because the condition $f^N\in (x^p, y^p) =\mathfrak m^{[p]}$ is clearly invariant under any linear change of  coordinates on $k[x, y]$, including the action of swapping $x$ and $y$; such a swap induces the action on  $k[a_0, \dots, a_d]$ that interchanges $a_i$ and $a_{d-i}$ for all $i$ and swaps 
$\mathcal C_{N, j}$ with  $\mathcal C_{N\, dN-j}$ for all $j$. This says that the greatest common square-free divisor $\mathcal C$ is invariant under swapping $a_i$ and $a_{d-i}$ for all $i$.

To prove that the closed set $X_{\leq \frac{N}{p}}$  has codimension two or more in $\PSym$, it suffices to show that we can find a projective linear space $\Lambda$ in $\PSym$ such that 
 $$
X_{\leq \frac{N}{p}} \, \cap \,  \Lambda 
$$
is codimension two or more in $\Lambda$. For this, define,
  for any fixed index $j$,  the ideal
\begin{equation}\label{j}
\mathcal A_j =  \underbrace{(a_0, \dots, a_{j-1}, \widehat{a_j}, \widehat{a_{j+1}}, \widehat{a_{j+2}}, a_{j+3}, \dots, a_d).}_{\text{all variables $a_t$ in $k[a_0, \dots, a_d]$  except $a_j, a_{j+1}, a_{j+2}$}}
\end{equation}
Letting $
\Lambda_j = \mathbb V(\mathcal A_j) \subseteq \PSym$, we will show 
we can choose $j$ (depending on whether $d$ is even or odd), so that $X_{\leq \frac{N}{p}} \cap \Lambda_j$ has codimension two or more in $\Lambda_j\cong \mathbb P^2$.

\noindent
{\bf The case when $d$ is even:}
Let $j=\frac{d-2}{2}.$
We claim that
\begin{equation}\label{EvenClaim}
X_{\leq \frac{N}{p}} \, \cap \, \mathbb V(\mathcal A_j)\, \subseteq \ \mathbb V(a_j, a_{j+1}) \ \cup \ \mathbb V(a_{j+1}, a_{j+2}),
\end{equation}
which will complete the proof in the case where $d$ is even because it implies each component of $X_{\leq \frac{N}{p}} \, \cap \, \mathbb V(\mathcal A_j)$ is contained in either $ \mathbb V(a_j, a_{j+1}) $ or $ \mathbb V(a_{j+1}, a_{j+2}),$ and so has codimension two or more in $\Lambda_j$.
To prove (\ref{EvenClaim}), observe it is equivalent to  stating that if  
$$ f=x^{d-2-j}y^{j}g = (xy)^{\frac{d-2}{2}} g =  (xy)^{\frac{d-2}{2}}(\underbrace{a_j x^2+a_{j+1} xy + a_{j+2}y^2}_{g})
$$
satisfies $\fpt(f)\leq \frac{N}{p}$, then either $a_j=a_{j+1}=0$ or $a_{j+1}=a_{j+2}=0$. 
Note that $\frac{d}{2}< \frac{p}{N}$ (since reciprocally, $\frac{2}{d}>\frac{N}{p}$), so we can 
apply Lemma \ref{genL1} with $i=j=\frac{d-2}{2}$ to see immediately that $a_j=0$ or $a_{j+2}=0$.  Without loss of  generality, say that $a_j=0$. In this case, we can write 
$$ f= x^{\frac{d-2}{2}}y^{\frac{d}{2}} ({a_{j+1} x + a_{j+2}y}),
$$
and we can again apply Lemma \ref{genL1}, this time with $(i, j) = (\frac{d-2}{2}, \frac{d}{2})$  to see that either $a_{j+1}=0$, in which case (\ref{EvenClaim}) follows, or $a_{j+2}=0$. But in this latter case, we have $a_j=a_{j+2}=0$, so that $f=a_{j+1}(xy)^{\frac{d}{2}}$. But then $\fpt f= \frac{2}{d}>\frac{N}{p}$. This establishes  (\ref{EvenClaim}), and hence the proof of Theorem \ref{codim2} for even $d$.

\noindent
{\bf The case when $d$ is odd:}
Let $j=\frac{d-3}{2}.$ 
We first claim that
\begin{equation}\label{OddClaim}
X_{\leq \frac{N}{p}} \, \cap \, \mathbb V(\mathcal A_j)\, \subseteq \ \mathbb V(a_{j+2}) \ \cup \ \mathbb V(a_j, a_{j+1}),
\end{equation}
or  equivalently, that all polynomials of the form 
\begin{equation}\label{oddd}
 f=x^{d-2-j}{y^j}g = x^{\frac{d-1}{2}}y^{\frac{d-3}{2}}(\underbrace{a_j x^2+a_{j+1} xy + a_{j+2}y^2}_{g})
\end{equation}
with $\fpt(f)\leq \frac{N}{p}$  must satisfy either $a_{j+2}=0$ or both $a_j=a_{j+1}=0$. 

We can  apply  Lemma \ref{genL1} to the indices  $(i, j) = (\frac{d-1}{2}, \frac{d-3}{2})$ to see that if $\fpt(f)\leq \frac{N}{p}$, then $g$ must be divisible by $x$ or  $y$. This says that in the expression (\ref{oddd}), 
  either $a_{j+2}=0$, in which case (\ref{OddClaim}) follows, or $a_j=0$. In this latter case, when $a_j=0$, we can write 
$$
f =(xy)^{\frac{d-1}{2}}(a_{j+1} x + a_{j+2}y)
$$
and again apply Lemma  \ref{genL1}, this time with  $i= j =\frac{d-1}{2}$ and $g=a_{j+1} x + a_{j+2}y$. Again, we see that because $\fpt f\leq \frac{N}{p}$,  either $a_{j+1}$ or $a_{j+2}$ must be zero. In either case, (\ref{OddClaim}) follows.

Now, having established (\ref{OddClaim}), we interpret it as saying that  each codimension one component of $X_{\leq \frac{N}{p}} \, \cap \, \mathbb V(\mathcal A_j)$ must be contained in $\mathbb V(a_{j+2})\cap \mathbb V(\mathcal A_j)$, 
which means that their union $\mathbb V(\mathcal C) \, \cap \, \mathbb V(\mathcal A_j)$ is contained in $\mathbb V(a_{j+2})\cap \mathbb V(\mathcal A_j).$
In particular, if $\overline{ \mathcal C}$ denotes the image of $\mathcal C$ in the quotient ring $k[a_0, \dots, a_d]/\mathcal A_j = k[a_j, a_{j+1}, a_{j+2}]$, we have 
$$
\mathbb V(\overline{ \mathcal C}) \subseteq \mathbb V(a_{j+2})
$$
in the copy of $\mathbb P^2$ defined by $\mathcal A_j$.  This says that $a_{j+2}\in \sqrt(\overline{ \mathcal C}),$  and so some power of $a_{j+2}$ ($=a_{\frac{d+1}{2}}$) appears with non-zero coefficient  in the  homogeneous polynomial 
$\mathcal C\in k[a_0, a_1, \dots, a_d]$. 
Now property (C) above implies that the same power of $a_{\frac{d-1}{2}}$ appears  in $\mathcal C$ with  the same non-zero coefficient. 
But this contradicts  property (B), since $a_{\frac{d-1}{2}}$ and $a_{\frac{d+1}{2}}$ have different degrees in the weighted grading.
This completes the proof of Theorem \ref{codim2} in the case where $d$ is odd. 
\end{proof}


\section{Low and Special degrees}\label{Smalld}

In this section, we apply our results to identify the precise lists of rational numbers that occur as the $F$-pure threshold of a reduced polynomial  in two variables over an algebraically closed field in certain cases, 
including the cases of fixed small degree $d$ (in terms of the  congruence class of $p$ modulo $d$) as well as certain fixed values of $d$ relative to $p$.

A list of potential values for the F-pure threshold was described in  conditions (I), (II) and (III) in \S 5 
(following  \cite[Thm 4.4]{HNWZ}).  While we saw in  Example \ref{NotSufficient} that these conditions are {\it  not sufficient} to identify the rational numbers that are $F$-pure thresholds of reduced polynomials in $k[x,y]$ in general, 
 we verify that they are sufficient  for every characteristic $p$ and  degree $d$ up to eight. These computations showcase the use of Theorems \ref{generic},  \ref{order} and \ref{p1} to compute $F$-pure threshold, as well as a few additional tricks.

\medskip

We first point out some consequences  for certain $d$ and $p$: 

\begin{prop} \label{1modd} If $2p\equiv 2\mod d$ and $d\geq 4$, then all reduced polynomials of degree $d$ over a field of characteristic $p$ have F-pure threshold either $\frac{2}{d}$ or its truncation
$$
\frac{2p-2}{dp} 
$$
at the first spot. {\bf Both} of these values occur as $F$-pure thresholds of reduced polynomials in every characteristic.
\end{prop}

\begin{proof}
By  Corollary  \ref{genericCOR},  the generic value of the $F$-pure threshold is $\frac{2}{d}$. Note our assumptions on $d$ imply that $p\neq 2$, and so $p\not | d$. Thus Theorem \ref{order} implies that
 there is a reduced polynomial whose $F$-pure threshold  is the first truncation $
\frac{2p-2}{dp}$. No other truncation can be the $F$-pure threshold of a reduced polynomial by condition (I). 
\end{proof}

 \begin{prop} Fix $p>3$.
   Consider reduced $f\in k[x, y]$ of degree $p+1$. Then $\fpt(f)$ is
   \begin{enumerate}
   \item $\frac{2}{p+1} $, or 
   \item $\frac{2p^2-2}{p^2(p+1)} $ (which is the truncation at the second spot)
    or \item 
    $\frac{1}{p}$ (which is the truncation at the first spot). 
    \end{enumerate}
    All these cases occur for reduced polynomials of degree $d$.
   \end{prop}
\begin{proof}
 The generic value of the $F$-pure threshold is $\frac{2}{p+1}$ by Theorem \ref{genericCOR},  as there is no natural number $e$ such that $2p^e\equiv 1\mod (p+1)$. 
Next, note that 
$$
  2p^2-2 = 2(p+1)(p-1) \equiv  0 \mod (p+1), 
  $$
  where as  $p\not\equiv 1 \mod (p+1)$, so 
 Theorem \ref{p1} guarantees that
 the second truncation of $\frac{2}{p+1}$,  $\frac{2p^2-2}{p^2(d+1)}, $ occurs as the $F$-pure threshold of some reduced polynomial.
  Finally, the first truncation $\frac{1}{p}$ occurs for reduced Frobenius forms, such at $x^{p+1}+y^{p+1}$ by \cite[1.1]{extremal}.
  \end{proof}

We highlight one more special case, which follows from restriction  (III) in Section 5:

\begin{prop}\label{d=p} Let $k$ be a field of characteristic $p>0$. Let $f$ be a reduced homogeneous polynomial of degree $d$ where $d=p^t$ or $d=2p^t$ for some $t\in \mathbb N$.
Then 
$$
\fpt(f) = \frac{2}{d}.
$$
\end{prop}

The next several propositions show that the conditions (I), (II), (III) precisely describe the possible $F$-pure thresholds in degrees up to 8.  
We first  dispose of the case $d=4$:

\begin{prop}\label{d=4} The precise list of all $F$-pure thresholds that can and do occur as the $F$-pure threshold of a reduced polynomial of degree $4$ is:
\begin{enumerate}
 \item[{\rm (a).}]  In  every characteristic $p$, the generic value of the $F$-pure threshold is $\frac{1}{2}$; 
\item[{\rm (b).}]  If $p=2$,  all reduced polynomials have $F$-pure threshold $\frac{1}{2}$; 
\item[{\rm (c).}]  If $p\neq 2$, then there are also reduced polynomials with $F$-pure threshold $\langle\frac{1}{2}\rangle_1=\frac{p-1}{2p}$.
\end{enumerate}
\end{prop}

\begin{proof} By \cite[4.4]{HNWZ} (conditions (I), (II), and (III)), every $F$-pure threshold of a reduced 
degree four polynomial is one of the types listed above. To see that each actually occurs for a reduced polynomial, 
 note that  (b) follows from Proposition \ref{d=p}, and (a) and (c),  follow from Corollary \ref{genericCOR} and Theorem \ref{order}, respectively, as $2p^e\neq 1 \mod 4$ for any natural number  $e$ while $p^1\equiv 1 \mod 2$ for all odd $p$. 
\end{proof}

For degree five, a list of possible values of the F-pure threshold  were described by Hara and  Monsky (when $p\neq 5$); see \cite[3.9]{HaraMonsky} or  \cite[4.4]{HNWZ}.
 In fact, for every $p$, all these values {\it do}  occur as the $F$-pure threshold for a reduced polynomial of degree $5$:
\begin{prop}
The complete list of  values that can and do occur as the $\mathrm{fpt}(f)$ of a  reduced form in $k[x, y]$ of degree 5 is:
\begin{enumerate}\label{d=5}
\item[{\rm (a).}]   If $p=5$: \,\,  $\displaystyle \frac{2}{5}$;
\item[{\rm (b).}] If $p \equiv 1 \mod 5$: \,  $\displaystyle \frac{2}{5}$ and \ $\displaystyle \frac{2p-2}{5p}$;
\item[{\rm (c).}] If $p \equiv 2 \mod 5$: \, $\displaystyle \frac{2p^2-3}{5p^2}$ and \  $\displaystyle \frac{2p^3-1}{5p^3}$ for $p\geq 7$, but only the value $\frac{1}{4}$ when $p=2$;
\item[{\rm (d).}] If $p \equiv 3 \mod 5$: \, $\displaystyle \frac{2p-1}{5p}$;
\item[{\rm (e).}]If $p \equiv 4 \mod 5$: \,  $\displaystyle \frac{2}{5}$, $\displaystyle \frac{2p-3}{5}$, and \ $\displaystyle \frac{2p^2-2}{5p^2}$.

\end{enumerate}
\end{prop}

\begin{remark}
It is shown in \cite[4.8]{HNWZ}, that for each congruence class in Proposition \ref{d=5} and each formula listed above, there {\it exist} $p$ and reduced $f$ such that $\fpt(f)$ has the given formula. Proposition \ref{d=5} is  stronger: it ensures that this holds {\it for every} prime $p$.
\end{remark}

\begin{proof}
The existence of reduced polynomials for (a) is immediate  from  Proposition \ref{d=p}  and (b)  from Proposition \ref{1modd}.

For (c),  the generic value of the $F$-pure threshold is the third truncation, which is $\frac{1}{4}$ when $p=2$ and  $\frac{2p^3-1}{5p^3}$ otherwise.  
The first truncation does not occur by condition (III). The second truncation is again $\frac{1}{4}$ if $p=2$ and 
 $\langle \frac{2}{5}\rangle_2=\frac{2p^2-3}{5p^2}$ for larger $p$. This latter value occurs for a reduced polynomial because it is the F-pure threshold of $f=x^5+y^5$. To see this, we observe that $(x^5+y^5)^N\in (x^{p^2}, y^{p^2})$ where $N=\frac{2p^2-3}{5}$. Indeed, 
  the only monomials of the expansion of  $(x^5+y^5)^N$ that are not {\it a priori} in $(x^{p^2}, y^{p^2})$ are $x^{p^2-1}y^{p^2-2}$ and $x^{p^2-2}y^{p^2-1}$, and neither of these has both exponents divisible by $5$. (Alternatively, one can compute that $\fpt(x^5+y^5) = \frac{2p^2-3}{5p^2}$  using
  \cite[3.1]{hernandez.diag_hypersurf}.)

Statement (d) follows from  Corollary \ref{genericCOR}. 
In case (e), the values $\frac{2}{5}$ and its second truncation are $F$-pure thresholds of reduced polynomials by 
Corollary \ref{genericCOR} and Theorem \ref{p1}, respectively.
The first truncation, $\frac{2p-3}{5}$, is the $F$-pure threshold of $x^5+y^5$, arguing as in case (c) above or  using \cite[3.1]{hernandez.diag_hypersurf}.  
\end{proof}

\begin{prop}\label{d=6} The precise list of all rational numbers that can and do occur as the $F$-pure threshold of a reduced homogeneous 
polynomial in two variables of degree $6$ is:
\begin{enumerate}
\item[{\rm (a).}]  If $p=3$: \ $\displaystyle \frac{1}{3}$;
 \item[{\rm (b).}]
 If $p \equiv 1 \mod 3$: \, $\displaystyle \frac{1}{3}$ and $\displaystyle \frac{p-1}{3p}$; 
 \item[{\rm (c).}]    If  $p \equiv 2  \mod 3$: \, $\displaystyle \frac{1}{3}$, $\displaystyle \frac{p-2}{3p}$, and $\displaystyle \frac{p^2-1}{3p^2}$ for $p\geq 5$,  but only $\frac{1}{3}$ and $\frac{1}{4}$ if $p=2$.
\end{enumerate}
\end{prop}
\begin{proof} The list above is created by the necessary conditions (I), (II), (III), as verified in \cite[4.6]{HNWZ}. It remains to show that every listed formula is achieved for every valid $p$. 
Case (a)  follows from Proposition \ref{d=p}   and (b) from Proposition \ref{1modd}.  For case (c), the values $\frac{1}{3}$ and its second truncation are $F$-pure thresholds of reduced polynomials by  Corollary \ref{genericCOR} and Theorem \ref{p1}, respectively.
 The first truncation is zero when $p=2$. For larger $p$, the first truncation  is the $F$-pure threshold of a reduced polynomial of the form $f=x^6+ax^3y^3+y^6$. To show this, it suffices to show that there is a value of $a$, not $a=\pm 2$, such that
 $f^N\in (x^p, y^p)$ where $N=\frac{p-2}{3}$. Expanding out $f^N$ as a sum of monomials $x^\alpha y^\beta$ of degree $6N=2p-4$, and using the fact that both $\alpha$ and $\beta$ are multiples of $3$, we see that
 $$
 f^N \equiv  C(a)(xy)^{p-2}  \mod (x^p, y^p),
 $$
 where $C(a)$ is a monic polynomial of degree $N$ in $a$. So  $\fpt(f)\leq \frac{N}{p}$ if and only if $C(a)=0$. Note finally that $\pm 2$ are {\it not} roots of $C,$ since 
 $$[(x^3\pm y^3)^2]^N \equiv \binom{2N}{N} (xy)^N =  \binom{2N}{N} (xy)^{p-2} \not\in (x^p, y^p)$$ as $2N<p$.
 So there is at least value of $a$ for which $f$ is reduced and  $C(a)=0$. For such choice of $a$, $\fpt(f)=\frac{N}{p}$. 
\end{proof}

The $F$-pure threshold of a reduced polynomial of degree 7 is necessarily included on the list below, using
 the conditions (I), (II) and (III), and verified in \cite[4.7]{HNWZ}. The next proposition shows that, for {\it every} $p$,   (I), (II) and (III) are {\it sufficient} as well as necessary.

\begin{prop}\label{d=7}
For each prime $p>0$,  the complete list of numbers that can and do occur the $F$-pure threshold of a reduced homogeneous polynomial in $k[x,y]$ of degree 7:
\begin{enumerate}
\item[{\rm (a).}] 
 If $p= 7$: \ $\frac{2}{7}.$ 
\item[{\rm (b).}] If $p \equiv 1  \mod 7$: \,  $\displaystyle \frac{2}{7}$ and  $\displaystyle \frac{2p-2}{7p}$;

\item[{\rm (c).}]  If $p \equiv 2 \mod 7$: \, $\displaystyle \frac{2p-4}{7p}$ and $\displaystyle \frac{2p^2-1}{7p^2}$ for $p>2$, but only $\frac{1}{4}$ when $p=2$;

\item[{\rm (d).}] If $p \equiv 3   \mod 7$: \, $\displaystyle \frac{2p^2-4}{7p^2}$, $\displaystyle \frac{2p^3-5}{7p^3}$, and $\displaystyle \frac{2p^4-1}{7p^4}$ for $p>3$, \\ \, but only 
 $\frac{23}{81}$ and $\frac{2}{9}$ when $p=3$;

\item[{\rm (e).}]  If $p \equiv 4  \mod 7$: \,  $\displaystyle \frac{2p-1}{7p}$;

\item[{\rm (f).}]  If $p \equiv 5  \mod 7$: \, $\displaystyle \frac{2p-3}{7p}$ and $\displaystyle \frac{2p^2-1}{7p^2}$;

\item[{\rm (g).}]  If $p \equiv 6  \mod 7$: \,  $\displaystyle \frac{2}{7}$, $\displaystyle \frac{2p-5}{7p}$, and $\displaystyle \frac{2p^2-2}{7p^2}$.

\end{enumerate}
\end{prop}
\begin{proof}
We must justify that each listed value occurs as the $F$-pure threshold of a reduced polynomial in $k[x, y]$. 
Case  (a) and (b)  follow from Propositions \ref{d=p} and \ref{1modd}, respectively. 

For (c), the second truncation occurs  by Corollary \ref{genericCOR}. To see that the first truncation occurs, we can then 
argue as in the proof of Proposition \ref{d=6} (c)  to show there is a value of $a$, not $a=\pm 2$, such that the reduced polynomial $f=x(x^6+ax^3y^3+y^6)$  satisfies $f^N\in (x^p, y^p)$, where $N= \frac{2p-4}{7}$.
Indeed, expanding out $f^N$ as a sum of monomials $x^\alpha y^\beta$ and noting that $\beta$ must be a multiple of 3, we see that there is a unique pair $(\alpha, \beta) = (i, j)$ such that
$$
 f^N \equiv  C(a)x^iy^j \mod (x^p, y^p),
 $$
for some polynomial $C\in k[a]$. Since  $2N<p$, note $x^N(x^3\pm y^3))^{2N}\not\in (x^p, y^p)$, so that $\pm 2$ are {\it not} roots of $C(a)$. So choosing  $a$ to be any root of  $C$,  the reduced polynomial $f=x(x^6+ax^3y^3+y^6)$ satisfies $\fpt(f)=\frac{N}{p}= \frac{2p-4}{7p}.$

For (d), the fourth truncation is the generic $F$-pure threshold by Corollary \ref{genericCOR}. 
When $p=3$,  the base $p$ expansion of $\frac{2}{7} =\frac{0}{3}+\frac{2}{9}+\frac{0}{27}+\dots$, so the  only possible additional value to check  is  whether or not the truncation $\frac{2}{9}$ is an $F$-pure threshold; it is,
 because the polynomial $f=x(x^6+ax^3y^3+y^6)$   satisfies $f^2\in (x^9, y^9)$ when $a^2=-2$.
When $p>3$,  the third truncation is the $F$-pure threshold of $x^7+y^7$ by Hernandez's formula \cite[Thm 3.1]{hernandez.diag_hypersurf}.

In  case (d) when $p>3$, we claim there is reduced  polynomial  $f$ of the form  $x(x^6+ax^3y^3+y^6) $ such that $f^N\in (x^{p^2}, y^{p^2})$ for $N=\frac{2p^2-4}{7}$. Indeed,  only three monomials of degree $7N$,
 $$x^{p^2-1}y^{p^2-3}, \,\,x^{p^2-2}y^{p^2-2},\,\, x^{p^2-3}y^{p^2-1},$$
 fail to be in the ideal $(x^{p^2}, y^{p^2})$ but only the third has $y$  exponent divisible by $3$.
 So we must show that the coefficient $C(a)\in k[a]$  of the term $x^{p^2-3}y^{p^2-1}$ in the expansion of  $f^N$ has some root other than $\pm 2$. 
Note that $C(a)$  is constructed by summing over all choices of one monomial from each of the $N$-copies of $(x^6+ax^3y^3+y^6)$ so that the exponents on $y$ sum to $p^2-1$; in particular, the coefficient of $a^T$ in $C(a)$ is
 \begin{equation}\label{Coefficient}
\underbrace{
 \binom{N}{T}}_{  \#  \ ax^3y^3 \text{ choices}} \
 \underbrace
 {
 \binom{N-T}
{\frac{p^2-3T-1}{6}},}_{  \# \ y^6  \text{ choices } }
\end{equation} 
or  zero if ${\frac{p^2-3T-1}{6}}\not\in \mathbb N$. Now suppose
$T=\frac{5p^2-17}{21}$. Note first that 
$$
\frac{1}{6}(p^2-3T-1) = \frac{1}{6}(p^2-\frac{5p^2-17}{7}-1)=\frac{1}{6}(\frac{2p^2+10}{7})=\frac{1}{3}(\frac{p^2+5}{7})\in \mathbb N
$$
(since $p^2=1\mod 3$), so the term $a^T$ appears in $C(a)$ with coefficient (\ref{Coefficient}). To compute the coefficient of $a^T$,
 note that
$$
 N-T=\frac{2p^2-4}{7}-\frac{5p^2-17}{21} = \frac{p^2+5}{21} = \frac{1}{6}(p^2-3T-1),$$
  so the factor $ \binom{N-T}
{\frac{p^2-3T-1}{6}}$ in (\ref{Coefficient}) is $1$. For the other factor in (\ref{Coefficient}), use the base $p$ expansions 
$$
N=\frac{2p-6}{7}\ p + \frac{6p-4}{7} \,\,\,\,\,\,{\text{and}}\,\,\,\,\,\, N-T=\begin{cases}
\frac{p -10 }{21}\ p + \frac{10p+5}{21} & {\text{ if }} p\equiv 1 \mod 3\\
\frac{p -17 }{21}\ p + \frac{17 p+5}{21} & {\text{ if }} p\equiv 2 \mod 3,
\end{cases}
$$
and   Lucas's theorem{\footnote{which in this case says $\binom{\alpha  p + \beta}{\gamma p + \delta} = \binom{\alpha}{\gamma} \binom{\beta}{\delta}\mod p$, when $0\leq \alpha, \beta, \gamma, \delta\leq p-1$}}
to check that $\binom{N}{N-T}\neq 0$.  So the coefficient $\binom{N}{N-T}$  of $a^T$ in $C(a)$ is non-zero, and $C(a)$ must have roots.
To check that $\pm 2$ are not among the roots of $C$,  observe that
 $$x^N(x^3\pm y^3)^{2N}\equiv \pm\binom{2N}{\frac{p^2-1}{3}} x^{p^2-3}y^{p^2-1} \mod 
  (x^{p^2}, y^{p^2}),$$
and the binomial coefficient   $\binom{2N}{\frac{p^2-1}{3}}\neq 0$.
Indeed,  using the base $p$ expansions 
 $$2N= \frac{4p-5}{7} \ p + \frac{5p-8}{7} \,\,\,\, \,\,\,{\text{and}}\,\,\,\,\,\,\,\,
 \frac{p^2-1}{3} =
 \begin{cases} \frac{p-1}{3} \ p + \frac{p-1}{3} \,\, {\text{ if }} p\equiv 1\mod 3  \\
\frac{p-2}{3} \, p + \frac{2p-1}{3} \,\, {\text{ if }} p\equiv 2\mod 3,  \\
\end{cases}
$$
 Lucas's theorem implies  $\binom{2N}{\frac{p^2-1}{3}}\neq 0$.

For (e), the generic $F$-pure threshold is the first truncation, so this is both the maximal and minimal possible $F$-pure threshold.
For (f), the generic $F$-pure threshold is the second truncation by Corollary \ref{genericCOR}, while the first truncation is $\fpt(x^7+y^7)$,  by the same method as in Proposition \ref{d=5} (c). 
And for
 (g), the values  $\frac{2}{7}$ and the second truncation are achieved by 
 Corollary \ref{genericCOR} and Theorem \ref{p1}, respectively, while the first truncation is $\fpt(x^7+y^7)$,   as in Proposition \ref{d=5} (c). 
\end{proof}

\begin{prop}\label{d=8}
Fix a prime $p$. The complete list of rational numbers that are $F$-pure thresholds of a reduced homogeneous polynomial of degree eight in two variables  is
\begin{enumerate}
\item [{\rm (a).}] If $p=2$: \, $\displaystyle \frac{1}{4}$;
\item [{\rm (b).}] If $p \equiv 1$ (mod 4): $\frac{1}{4}$ and $\frac{p-1}{4p}$
\item 
[{\rm (c).}] If $p \equiv 3$ (mod 4): $\frac{1}{4}$, $\frac{p-3}{4p}$, and $\frac{p^2-1}{4p^2}$ for $p>3$, {but only $\frac{1}{4}$ and $\frac{2}{9}$ when $p=3$.}
\end{enumerate} 
\end{prop}

\begin{proof} The list of potential $F$-pure thresholds is easily constructed from conditions (I), (II) and (III). We must verify that each occurs (for every valid $p$).
Statement (a) follows from Proposition \ref{d=p}  and (b) from  Corollary \ref{1modd}. 

For (c), note that  $\frac{1}{4}$ and the second truncation  are the $F$-pure thresholds of reduced polynomials by  Corollary \ref{genericCOR} and  Theorem \ref{p1}, respectively. Note that when $p=3$, the first truncation of $\frac{1}{4}$ base 3 is zero, 
so we have found all possible $F$-pure thresholds. 
 For $p\neq 3$, 
it remains to show that there exists a reduced degree $8$  polynomial  $f$  such that $f^N\in (x^p, y^p),$ where $N= \frac{p-3}{4}$. 
For $p\equiv 7\mod 8$, one checks that $\fpt(x^8+y^8)$ is such a polynomial,  using the method of Proposition \ref{d=5} (c).
For  $p\equiv 3 \mod 8$,  one checks there is a reduced polynomial of the form  $f=xy(x^6+ax^3y^3+y^6)$ with $\fpt(f)=\frac{p-3}{4p}$ similarly to the proof of Proposition \ref{d=6} (c). 
  \end{proof}

\section{Lower Bounds on the F-pure Threshold}

To conclude, we discuss lower bounds on the $F$-pure threshold. First we note that the first non-zero truncation of the base $p$ expansion of $\frac{2}{d}$ always gives a lower bound:

\begin{prop}\label{LB} Fix a field $k$ of characteristic $p>0$. Let $f$ be a homogeneous reduced polynomial of degree $d\geq 2$ in $k[x_1, \dots, x_n]$, $n\geq 2$.
Then
\begin{equation}\label{BoundBelow}
\fpt(f)\geq \left[\frac{2}{d}\right]_e,  {\text{ where $e\in \mathbb N$ is minimal such that }}  \left[\frac{2}{d}\right]_e\neq 0.
\end{equation}
Here, $\left[\frac{2}{d}\right]_e$ denotes the $e$-th truncation of the {\bf possibly-terminating} base $p$ expansion of  $\frac{2}{d}.$
\end{prop}

 Here, by "possibly-terminating" base $p$ expansion of a positive rational number  $\lambda$ we mean that, if $\lambda$ has more than one base $p$ expansion, we choose the one that terminates. Writing $\lambda$ in lowest terms, note that unless the denominator is a power of $p$, then $\lambda$ has only one base $p$ expansion so 
 the "possibly-terminating" base $p$ expansion is simply its non-terminating expansion. 
  See Subsection \ref{TruncationBasics}.

\begin{proof}
The $F$-pure threshold can only decrease after modding out  linear forms \cite[3.6]{extremal}.
 Modulo $n-2$ general enough forms, we get a reduced polynomial  $\overline f$ in two variables whose $F$-pure threshold is a lower bound for $\fpt(f)$. If $d=p^e$ or $2p^e$, the $F$-pure threshold of $\overline f$ is $\frac{2}{d}$ by Proposition \ref{d=p}, which has terminating base   $p$ expansion $\frac{2}{p^e}$ or $\frac{1}{p^e}$, respectively, showing the lower bound is satisfied. Otherwise, $\frac{2}{d}$ does not have denominator a power of $p$, so its possibly-terminating and its non-terminating base $p$ expansion are the same. Since we know that the $F$-pure threshold of the reduced polynomial $\overline f$ 
  in two variables is  either $\frac{2}{d}$ or some 
 truncation of the non-terminating base $p$ expansion by \cite[4.4]{HNWZ}, the first non-zero truncation is certainly a lower bound.
\end{proof}

\begin{cor}\label{sharpLB}
The bound  of Proposition \ref{LB} is sharp for homogeneous polynomials of degree $d$ in the range $p^e+1\leq d \leq 2p^e$ in $k[x_1, \dots, x_n]$ ($n\geq 2$). That is, there are reduced polynomials of degree $d$ whose $F$-pure threshold is $ \left[\frac{2}{d}\right]_e$.
\end{cor}

\begin{proof} For $d$ in the stated range, it is easy to compute that $ \left[\frac{2}{d}\right]_e =\frac{1}{p^e}$ and that all previous truncations are zero. Thus we only need to find reduced polynomials of degree  $d$  with $F$-pure threshold $\frac{1}{p^e}$. Since we already know this is a lower bound (Proposition \ref{LB}), we need only show that there are are reduced polynomials $f$  with $\fpt(f)\leq \frac{1}{p^e}$. But 
any reduced  polynomial $f$ in $(x_1^{p^e}, \dots,  x_n^{p^e})$ satisfies 
\begin{equation}\label{1}
\fpt(f) \leq \frac{1}{p^e},
\end{equation} by definition of $F$-pure threshold. So for such $f$, the lower bound $ \left[\frac{2}{d}\right]_e$ of Proposition \ref{LB} is sharp.
\end{proof}

Putting together results proved in earlier sections, we summarizes cases where we know the lower bound is sharp:
\begin{cor}\label{LBcor}
The bound of Proposition \ref{LB} is sharp for reduced polynomials of degree $d\geq 4$ in $k[x_1, \dots, x_n]$ in each of the following cases:
\begin{enumerate}
\item $p^e+1\leq d \leq 2p^e$;
\item $2p=1\mod d$;
\item $2p=2\mod d$;
\item $d\leq 8$.
\end{enumerate}
\end{cor}

\begin{proof}
Simply note that a reduced polynomial in two variables can be viewed as a reduced polynomial in $k[x_1, \dots, x_n]$, and apply  Corollary \ref{sharpLB}, Proposition \ref{1modd}, Theorem \ref{p1} respectively, as well as the computations for degrees $4\leq d\leq 8$ in Section 7.
\end{proof}

\begin{remark}
The lower bound is not sharp in general: in the case where $d=3$ and $p\equiv 1\mod p$, no reduced polynomial has $F$-pure threshold {\it any} truncation of $\frac{2}{3}$.
\end{remark}

In all cases we know where the lower bound of  Proposition \ref{LB} fails, the reason is that the first non-zero truncation fails one of the conditions (I), (II) or (III) from Section 5 (that is, from \cite[4.4]{HNWZ}). 
We suspect  that the {\it first    truncation of the possibly-terminating  base $p$ expansion of $\frac{2}{d}$ which satisfies (I) , (II) and (III) is {\bf{ precisely}} the minimal possible $F$-pure threshold of a reduced polynomial of degree $d\geq 2$
in $k[x_1, \dots, x_n]$ in general.}

We do not know whether  the first
  truncation of the possibly-terminating  base $p$ expansion of $\min(1, \frac{n}{d})$ which satisfies (I), (II) and (III) could be the minimal possible $F$-pure threshold for  polynomials in $k[x_1, \dots, x_n]$ defining a smooth hypersurface in $\mathbb P^{n-1}.$
  
  \begin{remark}
  In \cite{extremal}, a lower bound of $\frac{1}{d-1}$ was proved for a reduced homogeneous polynomial of degree $d$ (in any number of variables); moreover, this bound was shown to be sharp only when $d=p^e+1$ for some $e\in \mathbb N$.
  This comports with Corollary \ref{LBcor} (1) above.
  
  \end{remark} 

  \begin{remark} A lower bound on the $F$-pure threshold of elements $f$ in more  general local rings $(R, \m)$ is proved by N\'u\~nez-Betancourt and Smirnov: $\fpt(f)\geq \frac{1}{e_{HK}(R/fR)}$ where $e_{HK}(R/fR)$ is the Hilbert-Kunz  multiplicity of
  the quotient local ring $R/fR$  \cite[2.2]{Luis.Ilya}. In our main case of interest, however, this  is much weaker than the bound  in Proposition \ref{LB}: when $R=k[x,y]$,  the Hilbert-Kunz multiplicity of the one dimensional ring $R/(f)$ is the same as its Hilbert-Samuel multiplicity (see {\it e.g} \cite[3.3]{Huneke}), so the  N\'u\~nez-Betancourt-Smirnov bound simply recovers  Proposition \ref{basic} (4).
  \end{remark}

\bibliographystyle{amsalpha}
\bibliography{bibdatabase.bib}
\end{document}